\documentclass[a4paper,11pt]{article}
\usepackage[latin1]{inputenc}
\usepackage[T1]{fontenc}
\usepackage[english]{babel}
\usepackage{amssymb}
\usepackage{amsfonts, amsmath}
\newtheorem{Theo}{Theorem}

\newtheorem{Rem}{Remark}[section]
\newtheorem{Exmp}{Example}[section]
\title{Integrable systems and complex geometry}
\author{\textbf{A. Lesfari}
\\\emph{Department of Mathematics}
\\\emph{Faculty of Sciences}
\\\emph{University of Choua\"{i}b Doukkali}
\\\emph{B.P. 20, El-Jadida, Morocco}.
\\\emph{E. mail : lesfariahmed@yahoo.fr, lesfari@ucd.ac.ma}}
\date{}
\begin{document}
\maketitle
\begin{abstract}
In this paper, we discuss an interaction between complex geometry
and integrable systems. Section 1 reviews the classical results on
integrable systems. New examples of integrable systems, which have
been discovered, are based on the Lax representation of the
equations of motion. These systems can be realized as straight
line motions on a Jacobi variety of a so-called spectral curve. In
section 2, we study a Lie algebra theoretical method leading to
integrable systems and we apply the method to several problems. In
section 3, we discuss the concept of the algebraic complete
integrability (a.c.i.) of hamiltonian systems. Algebraic
integrability means that the system is completely integrable in
the sens of the phase space being folited by tori, which in
addition are real parts of a complex algebraic tori (abelian
varieties). The method is devoted to illustrate how to decide
about the a.c.i. of hamiltonian systems and is applied to some
examples. Finally, in section 4 we study an a.c.i. in the
generalized sense which appears as covering of a.c.i. system. The
manifold invariant by the complex flow is
covering of abelian variety.\\
\emph{Mathematics Subject Classification} (2000). 37J35, 70H06,
14H40, 14H70, 14M10.

\end{abstract}

\section{Integrable systems}

Let $M$ be an even-dimensional differentiable manifold. A
symplectic structure (or symplectic form) on $M$ is a closed
non-degenerate differential 2-form $\omega$ defined everywhere on
$M.$ The non-degeneracy condition means that
$$\forall x\in M,\forall \xi \neq 0,\exists \eta :
\omega \left( \xi ,\eta \right) \neq 0,\left( \xi ,\eta \in
T_{x}M\right) .$$ The pair $\left( M,\omega \right) $ is called a
symplectic manifold.
\begin{Exmp}
The cotangent bundle $T^{*}M$ possesses in a natural way a
symplectic structure. In a local coordinate $\left( x_{1},\ldots
,x_{n},y_{1},\ldots ,y_{n}\right),$ $2n=\dim M,$ the form $\omega$
is given by $\omega =\sum_{k=1}^{n}dx_{k}\wedge dy_{k}.$
\end{Exmp}
\begin{Exmp}
Another important class of symplectic manifolds consists of the
coadjoints orbits $\mathcal{O}\subset \mathcal{G}^*,$ where
$\mathcal{G}$ is the algebra of a Lie group $\mathcal{G}$ and
$\mathcal{G}_\mu=\{Ad^*_g \mu :g\in \mathcal{G}\}$ is the orbit of
$\mu \in \mathcal{G}^*$ under the coadjoint representation.
\end{Exmp}

\begin{Theo}
a) Let $I:T_{x}^{*}M\longrightarrow T_{x}M,\text{ }\omega _{\xi
}^{1}\longmapsto \xi ,$ be a map defined by $\omega _{\xi
}^{1}\left( \eta \right) =\omega \left( \eta ,\xi \right) ,\text{
}\forall \eta \in T_{x}M.$ Then $I$ is an isomorphism generated by
the symplectic form $\omega.$\\
b) The symplectic form $\omega $ induces a hamiltonian vector
field $IdH:M\longrightarrow T_{x}M,\text{ }x\longmapsto IdH\left(
x\right),$ where $H:M\longrightarrow \mathbb{R},$ is a
differentiable function (called hamiltonian). In others words, the
differential system defined by
$$\dot x(t)=X_{H}\left( x\left( t\right) \right)=
IdH\left( x\right),$$ is a hamiltonian vector field associated  to
the function $H.$ The matrix that is associated to an hamiltonian
system determine a symplectic structure.
\end{Theo}
\emph{Proof}. a) Denote by $I^{-1}$ the map
$I^{-1}:T_{x}M\longrightarrow T_{x}^{*}M, \text{ }\xi \longmapsto
I^{-1}\left( \xi \right) \equiv \omega _{\xi}^{1},$ with
$I^{-1}\left( \xi \right) \left( \eta \right) =\omega _{\xi
}^{1}\left( \eta \right) =\omega \left( \eta ,\xi \right) ,\text{
}\forall \eta \in T_{x}M.$ The fact that the form $\omega $ is
bilinear implies that
\begin{eqnarray}
I^{-1}\left( \xi _{1}+\xi _{2}\right) \left( \eta \right)
&=&\omega \left( \eta ,\xi _{1}+\xi _{2}\right) ,\nonumber\\
&=&\omega \left( \eta ,\xi _{1}\right) +\omega \left( \eta ,\xi
_{2}\right) ,\nonumber\\
&=&I^{-1}\left( \xi _{1}\right) \left( \eta \right) +I^{-1}\left(
\xi _{2}\right) \left( \eta \right) ,\text{ }\forall \eta \in
T_{x}M.\nonumber
\end{eqnarray}
Now, since $\dim T_{x}M=\dim T_{x}^{*}M,$ to show that $I^{-1}$ is
bijective, it suffices to show that is injective.  The form
$\omega$ is non-degenerate, it follows that
$$
KerI^{-1}=\left\{ \xi \in T_{x}M:\omega \left( \eta ,\xi \right)
=0, \text{ }\forall \eta \in T_{x}M\right\} =\left\{ 0\right\}.$$
Hence $I^{-1}$ is an isomorphism and consequently $I$ is also an
isomorphism (the inverse of an isomorphism is an isomorphism).\\
b) Let $\left(x_{1},\ldots ,x_{m}\right) $ be a local coordinate
system on $M,$ $(m=\dim M).$ We have
\begin{equation}\label{eqn:euler}
\dot x(t)=\sum_{k=1}^{n}\frac{\partial H} {\partial x_{k}}I\left(
dx_{k}\right) =\sum_{k=1}^{n}\frac{\partial H}{\partial
x_{k}}\xi^{k},
\end{equation}
where $I\left( dx_{k}\right) =\xi ^{k}\in T_{x}M$ is defined such
that :$\forall \eta \in T_{x}M,\text{ }\eta _{k}=dx_{k}\left( \eta
\right)= \omega \left( \eta ,\xi ^{k}\right),$ (\text{k-th
component of}\text{}$\eta$). Define $\left( \eta _{1},\ldots ,\eta
_{m}\right) $ and $\left( \xi _{1}^{k},\ldots ,\xi _{m}^{k}\right)
$ to be respectively the components of $\eta $ and $\xi ^{k}$,
then
\begin{eqnarray}
\eta _{k}&=&\omega \left( \sum_{i=1}^{m}\eta _{i}\frac{\partial
}{\partial x_{i}},
\sum_{j=1}^{m}\xi _{j}^{k}\frac{\partial }{\partial x_{j}}\right) ,\nonumber\\
&=&\sum_{i=1}^{m}\eta _{i}\left( \frac{\partial }{\partial x_{i}},
\frac{\partial }{\partial x_{j}}\right) \xi _{j}^{k},\nonumber\\
&=&\left( \eta _{1},\ldots ,\eta _{m}\right) J^{-1}
\left(\begin{array}{c}
\xi _{1}^{k}\\
\vdots \\
\xi _{m}^{k}
\end{array}
\right) ,\nonumber
\end{eqnarray}
where $J^{-1}$ is the matrix defined by $J^{-1}\equiv \left(
\omega \left( \frac{\partial }{\partial x_{i}}, \frac{\partial
}{\partial x_{j}}\right) \right) _{1\leq i,j\leq m}.$ Since this
matrix is invertible\footnote{Indeed, it suffices to show that the
matrix $J^{-1}$ has maximal rank. Suppose this were not possible,
i.e., we assume that $rank(J^{-1})\neq m.$ Hence
$\sum_{i=1}^{m}a_{i}\omega \left( \frac{\partial }{\partial
x_{i}},\frac{\partial }{\partial x_{j}}\right) =0,$ $\forall 1\leq
j\leq m,$ with $a_{i}$ not all null and $\omega \left(
\sum_{i=1}^{m}a_{i}\frac{\partial }{\partial x_{i}},\frac{\partial
}{\partial x_{j}}\right) =0,$ $\forall 1\leq j\leq m.$ In fact,
since $\omega $ is non-degenerate, we have
$\sum_{i=1}^{m}a_{i}\frac{\partial }{\partial x_{i}}=0.$ Now
$\left( \frac{\partial }{\partial x_{1}},\ldots ,\frac{\partial
}{\partial x_{m}}\right) $ is a basis of $T_{x}M,$ then $a_{i}=0,$
$\forall i,$ contradiction.}, we can search $\xi ^{k}$ such that :
$$J^{-1}\left(\begin{array}{c}
\xi _{1}^{k}\\
\vdots \\
\xi _{m}^{k}
\end{array}
\right)= \left(\begin{array}{cc}
0&\\
\vdots \\
0&\\
1&\leftrightsquigarrow \text{k-th place}\\
0&\\
\vdots \\
0&
\end{array}
\right).$$ The matrix $J^{-1}$ is invertible, which implies
$$\left(\begin{array}{c}
\xi _{1}^{k}\\
\vdots \\
\xi _{m}^{k}
\end{array}
\right)= J\left(\begin{array}{c}
0\\
\vdots \\
0\\
1\\
0\\
\vdots \\
0
\end{array}
\right),$$ from which $\xi ^{k}=$(k-th column of J), i.e., $\xi
_{i}^{k}=J_{ik},$ $1\leq i\leq m,$ and consequently $\xi
^{k}=\sum_{i=1}^{m}J_{ik}\frac{\partial }{\partial x_{i}}.$ It is
easily verified that the matrix $J$ is
skew-symmetric\footnote{Indeed, since $\omega $ is symmetric i.e.,
$\omega \left( \frac{\partial }{\partial x_{i}},\frac{\partial
}{\partial x_{j}}\right) =-\omega \left( \frac{\partial }{\partial
x_{j}},\frac{\partial }{\partial x_{i}}\right) ,$ it follows that
$J^{-1}$ is skew-symmetric. Then, $I=J.J^{-1}=\left( J^{-1}\right)
^{\top }.J^{\top }=-J^{-1}.J,$ and consequently $J^{\top }=J.$}.
From $\left( 1\right) $ we deduce that
$$\dot x(t)=\sum_{k=1}^{m}\frac{\partial
H}{\partial x_{k}}\sum_{i=1}^{m}J_{ik}\frac{\partial }{\partial
x_{i}}=\sum_{i=1}^{m}\left( \sum_{k=1}^{m}J_{ik}\frac{\partial
H}{\partial x_{k}}\right) \frac{\partial }{\partial x_{i}}.$$
Writing $\dot x(t)=\sum_{i=1}^{m}\frac{dx_{i}\left( t\right)
}{dt}\frac{\partial }{\partial x_{i}},$ it is seen that $\dot
x_{i}\left( t\right)= \sum_{k=1}^{m}J_{ik}\frac{\partial
H}{\partial x_{k}},1\leq i\leq j\leq m,$ which can be written in
more compact form $\dot x(t)= J\left( x\right) \frac{\partial
H}{\partial x},$ this is the hamiltonian vector field associated
to the function $H.$ This
concludes the proof of the theorem.\\

We define a Poisson bracket (or Poisson structure) on the space
$\mathcal{C}^{\infty }$ as
$$\left\{ ,\right\} :\mathcal{C}^{\infty }\left( M\right)
\times \mathcal{C}^{\infty }\left( M\right) \longrightarrow
\mathcal{C}^{\infty }\left( M\right) ,\text{ }\left( F,G\right)
\longmapsto \left\{ F,G\right\} ,$$ where $\left\{ F,G\right\}
=d_{u}F\left( X_{G}\right) =X_{G}F\left( u\right) = \omega \left(
X_{G},X_{F}\right).$ This bracket is skew-symmetric $\left\{
F,G\right\} =-\left\{ G,F\right\} ,$ obeys the Leibniz rule
$\left\{ FG,H\right\} =F\left\{ G,H\right\} +G\left\{ F,H\right\}
,$ and satisfies the Jacobi identity
$$\left\{ \left\{ H,F\right\} ,G\right\} +\left\{ \left\{ F,G\right\} ,
H\right\} +\left\{ \left\{ G,H\right\} ,F\right\} =0.$$ When this
Poisson structure is non-degenerate, we obtain the
symplectic structure discussed above.\\

Consider now $M=\mathbb{R}^{n}\times \mathbb{R}^{n}$ and let $p\in
M.$ By Darboux's theorem $\left[3\right] $, there exists a local
coordinate system$\left(x_{1},\ldots ,x_{n}, y_{1},\ldots
,y_{n}\right) $ in a neighbourhood of $p$ such that $$\{H,F\}
=\sum_{i=1}^{n}\left( \frac{\partial H}{\partial
x_{i}}\frac{\partial F}{\partial y_{i}}-\frac{\partial H}{\partial
y_{i}}\frac{\partial F}{\partial x_{i}}\right).$$ Then
$X_{H}=\sum_{i=1}^{n}\left( \frac{\partial H}{\partial
x_{i}}\frac{\partial }{\partial y_{i}}-\frac{\partial H}{\partial
y_{i}}\frac{\partial }{\partial x_{i}}\right),$ and
$X_{H}F=\left\{ H,F\right\},\text{ }\forall F\in
\mathcal{C}^{\infty }\left( M\right).$ A nonconstant function $F$
is called an integral (first integral or constant of motion) of
$X_F,$ if $X_{H}F=0.$ In particular, $H$ is integral. Two
functions $F$ and $G$ are said to be in involution or to commute,
if $\left\{ F,G\right\}=0.$ The hamiltonian systems form a Lie
algebra.\\

We now give the following definition of the Poisson bracket :
$$\left\{ F,G\right\}=
\left\langle \frac{\partial F}{\partial x}, J\frac{\partial
G}{\partial x}\right\rangle= \sum_{i,j}J_{ij}\frac{\partial
F}{\partial x_{i}}\frac{\partial G}{\partial x_{j}}.$$ After some
algebraic manipulation, we deduce that If
$$\sum_{k=1}^{2n}\left( J_{kj}\frac{\partial J_{li}}
{\partial x_{k}}+J_{ki}\frac{\partial J_{jl}} {\partial
x_{k}}+J_{kl}\frac{\partial J_{ij}}{\partial
x_{k}}\right)=0,\text{ }\forall 1\leq i,j,l\leq 2n,$$ then $J$\
satisfies the Jacobi identity.\\

Consequently, we have a complete characterization of hamiltonian
vector field
\begin{equation}\label{eqn:euler}
\dot x(t)=X_{H}\left( x\left( t\right) \right) =J\frac{\partial
H}{\partial x},\text{ }x\in M,
\end{equation}
where $H:M\longrightarrow \mathbb{R},$ is a differentiable
function (the hamiltonian) and $J=J\left( x\right) $ is a
skew-symmetric matrix, possibly depending on $x\in M,$ for which
the corresponding Poisson bracket satisfies the Jacobi identity :
$$\left\{ \left\{ H,F\right\} ,G\right\} +\left\{ \left\{ F,G\right\}
,H\right\} +\left\{ \left\{ G,H\right\} ,F\right\}=0,$$ with
$\left\{ H,F\right\} =\left\langle \frac{\partial H}{\partial
x},J\frac{\partial F} {\partial x}\right\rangle
=\sum_{i,j}J_{ij}\frac{\partial H}{\partial x_{i}} \frac{\partial
F}{\partial x_{j}},$ the Poisson bracket.

\begin{Exmp}
An important special case is when $J=\left(\begin{array}{cc}
O&-I\\
I&O
\end{array}\right),$
where $I$ is the $n\times n$ identity matrix. The condition on $J$
is trivially satisfied. Indeed, here the matrix $J$ do not depend
on the variable $x$ and we have
$$
\left\{ H,F\right\} =\sum_{i=1}^{2n}\frac{\partial H}{\partial
x_{i}}\sum_{j=1}^{2n}J_{ij}\frac{\partial F}{\partial
x_{j}}=\sum_{i=1}^{n}\left(\frac{\partial H}{\partial
x_{n+i}}\frac{\partial F}{\partial x_{i}}-\frac{\partial
H}{\partial x_{i}}\frac{\partial F}{\partial x_{n+i}}\right) .$$
Moreover, equations (10) are transformed into
$$\dot q_{1}=\frac{\partial H}{\partial p_{1}},\ldots ,\dot q_{n}
=\frac{\partial H}{\partial p_{n}},\text{ }\dot
p_{1}=-\frac{\partial H} {\partial q_{1}},\ldots ,\dot
p_{n}=-\frac{\partial H}{\partial q_{n}},$$ o\`{u}
$q_{1}=x_{1},\ldots ,q_{n}=x_{n},p_{1}=x_{n+1},\ldots
,p_{n}=x_{2n}.$ These are exactly the well known differential
equations of classical mechanics in canonical form.
\end{Exmp}

It is a fundamental and important problem to investigate the
integrability of hamiltonian systems. Recently there has been much
effort given for finding integrable hamiltonian systems, not only
because they have been on the subject of powerful and beautiful
theories of mathematics, but also because the concepts of
integrability have been applied to an increasing number of applied
sciences. The so-called Arnold-Liouville theorem play a crucial
role in the study of such systems; the regular compact level
manifolds defined by the intersection of the constants of motion
are diffeomorphic to a real torus on which the motion is
quasi-periodic as a consequence of the following purely
differential geometric fact : a compact and connected
$n$-dimensional manifold on which there exist $n$ vector fields
which commute and are independent at every point is diffeomorphic
to an $n$-dimensional real torus and each vector field will define
a linear flow there.

\begin{Theo}
(Arnold-Liouville theorem)$\left[3,14\right] $\ : Let
$H_{1}=H,H_{2},...,H_{n},$ be $n$ first integrals on a
2n-dimensional symplectic manifold that are functionally
independent (i.e., $dH_1\wedge ... \wedge dH_n\neq 0$), and
pairwise in involution. For generic $c=(c_1,...,c_n)$ the level
set
$$M_{c}= \bigcap_{i=1}^{n}\left\{ x\in M:H_{i}\left( x\right)
=c_{i},\text{ }c_{i}\in \mathbb{R}\right\} ,$$ will be an
n-manifold. If $M_{c}$ is compact and connected, it is
diffeomorphic to an $n$-dimensional torus
$\mathbb{T}^n=\mathbb{R}^{n}/\mathbb{Z}^{n}$ and the solutions of
the system (2) are then straight-line motions on $\mathbb{T}^n$.
If $M_c$ is not compact but the flow of each of the vector fields
$X_{H_k}$ is complete on $M_c,$ then $M_c$ is diffeomorphic to a
cylinder $\mathbb{R}^k\times \mathbb{T}^{n-k}$ under which the
vector fields $X_{H_k}$ are mapped to linear vector fields.
\end{Theo}

As a consequence, we obtain the concept of complete integrability
of a hamiltonian system. For the sake of clarity, we shall
distinguish
two cases :\\
$\textbf{a})$ \underline{Case $1$} : $d\acute{e}t$ $J\neq 0.$ The
rank of the matrix $J$ is even, $m=2n.$ A hamiltonian system
$\left( 2\right) $ is completely integrable or
Liouville-integrable if there exist $n$ firsts integrals
$H_{1}=H,H_{2},\ldots ,H_{n}$ in involution, i.e.,
$\{H_{k},H_{l}\}=0,\text{ }1\leq k,l\leq n,$ with linearly
independent gradients, i.e., $dH_{1}\wedge ...\wedge dH_{n}\neq
0.$ For generic $c=(c_1,...,c_n)$ the level set
$$M_{c}= \bigcap_{i=1}^{n}\left\{ x\in M:H_{i}\left( x\right)
=c_{i},\text{ }c_{i}\in \mathbb{R}\right\} ,$$ will be an
n-manifold. By the Arnold-Liouville theorem, if $M_{c}$ is compact
and connected, it is diffeomorphic to an $n$-dimensional torus
$\mathbb{T}^n=\mathbb{R}^{n}/\mathbb{Z}^{n}$ and each vector field
will define a linear flow there. In some
open neighbourhood of the torus there are coordinates $%
s_{1},\ldots ,s_{n},\varphi _{1},\ldots ,\varphi _{n}$ in which
$\omega $ takes the form $ \omega =\sum_{k=1}^{n}ds_{k}\wedge
d\varphi _{k}.$ Here the functions $s_{k}$ (called
action-variables) give coordinates in the direction transverse to
the torus and can be expressed functionally in terms of the firsts
integrals $H_{k}.$ The functions $\varphi _{k}$ (called
angle-variables) give standard angular coordinates on the torus,
and every vector field $X_{H_{k}}$ can be written in the form
$\dot \varphi_{k}=h_{k}\left( s_{1},\ldots ,s_{n}\right) ,$ that
is, its integral trajectories define a conditionally-periodic
motion on the torus. In a neighbourhood of the torus the
hamiltonian vector field $ X_{H_{k}}$ take the following form $
\dot s_{k}=0, \dot \varphi_{k}=h_{k}\left( s_{1},\ldots
,s_{n}\right) ,$
and can be solved by quadratures.\\
$\textbf{b})$ \underline{Case $2$} : $d\acute{e}t$ $J=0.$ We
reduce the problem to $m=2n+k$ and we look for $k$ Casimir
functions (or trivial invariants) $H_{n+1},...,H_{n+k},$ leading
to identically zero hamiltonian vector fields $J\frac{\partial
H_{n+i}}{\partial x}=0,\text{ }1\leq i\leq k.$ In other words, the
system is hamiltonian on a generic symplectic manifold
$$\bigcap_{i=n+1}^{n+k}\left\{ x\in \mathbb{R}^{m}:H_{i}\left(
x\right) =c_{i}\right\},$$ of dimension $m-k=2n.$ If for most
values of $c_i\in \mathbb{R},$ the invariant manifolds
$$\bigcap_{i=1}^{n+k}\left\{ x\in \mathbb{R}^{m}:H_{i}\left( x\right)
=c_{i}\right\},$$ are compact and connected, then they are
n-dimensional tori $\mathbb{T}^n=\mathbb{R}^n/\mathbb{Z}^n$ by the
Arnold-Liouville theorem and the hamiltonian flow is linear in
angular coordinates of the torus.

\section{Isospectral deformation method}

A Lax equation is given by a differential equation of the form
\begin{equation}\label{eqn:euler}
\dot A\left( t\right) =\left[ A\left( t\right) , B\left( t\right)
\right] \text{ or }\left[ B\left( t\right) ,A\left( t\right)
\right],
\end{equation}
where
$$A\left( t\right) =\sum_{k=1}^{N}A_{k}\left( t\right) h^{k},
\text{ }B\left( t\right) =\sum_{k=1}^{N}B_{k}\left( t\right)
h^{k},$$ are functions depending on a parameter $h$ (spectral
parameter) whose coefficients $A_{k}$ and $B_{k}$ are matrices in
Lie algebras. The pair $\left( A,B\right) $ is called Lax pair.
This equation established a link between the Lie group theoretical
and the algebraic geometric approaches to complete integrability.
The solution to (3) has the form $A(t)=g(t)A(0)g(t)^{-1},$ where
$g(t)$ is a matrix defined as $\dot g\left( t\right) =-A(t)g(t).$
We form the polynomial $P\left( h,z\right) =\det \left(
A-zI\right),$ where $z$ is another variable and $I$ the $n\times
n$ identity matrix. We define the curve (spectral curve)
$\mathcal{C},$ to be the normalization of the complete algebraic
curve whose affine equation is $ P\left( h,z\right) =0.$

\begin{Theo}
The polynomial $P\left( h,z\right)$ is independent of $t.$
Moreover, the functions $tr\left( A^{n}\right) $ are first
integrals for (3).
\end{Theo}
\emph{Proof}. Let us call $L\equiv A-zI.$ Observe that
$$
\dot P=\det L.tr\left( L^{-1}\dot L\right) =\det L.tr\left(
L^{-1}BL-B\right)=0,$$ since $trL^{-1}BL=trB.$ On the other hand
\begin{eqnarray}
\dot  A^{n} &=&\dot A A^{n-1}+
A\dot A A^{n-2}+\cdots +A^{n-1}\dot A,\nonumber\\
&=&\left[A,B\right] A^{n-1}+A\left[A,B\right] A^{n-2}+\cdots
+A^{n-1}\left[ A,B\right],\nonumber\\
&=&\left( AB-BA\right) A^{n-1}+\cdots +A^{n-1}\left( AB-BA\right),\nonumber\\
&=&ABA^{n-1}-BA^{n}+\cdots +A^{n}B-A^{n-1}BA,\nonumber\\
&=&A\left( BA^{n-1}\right) -\left( BA^{n-1}\right) A+\cdots
+A\left( A^{n-1}B\right) -\left( A^{n-1}B\right) A.\nonumber
\end{eqnarray}
Since $tr\left( X+Y\right) =trX+trY,$ $trXY=trYX,$ $X,Y\in
\mathcal{M}_{n}\left( \mathbb{C}\right) ,$ we obtain
$$\frac{d}{dt}tr\left( A_{h}^{n}\right)=
tr\frac{d}{dt}\left( A_{h}^{n}\right) =0,$$ and consequently
$tr\left( A^{n}\right) $ are first integrals of motion. This
ends the proof of the theorem.\\

We have shown that a hamiltonian flow of the type $\left( 3\right)
$ preserves the spectrum of $A$ and therefore its characteristic
polynomial. The curve $\mathcal{C}:P(z,h)=\det \left(
A(h)-zI\right)=0,$ is time independent, i.e., its coefficients
$tr\left( A^{n}\right)$ are integrals of the motion (equivalently,
$A(t)$ undergoes an isospectral deformation. Some hamiltonian
flows on Kostant-Kirillov coadjoint orbits in subalgebras of
infinite dimensional Lie algebras (Kac-Moody Lie algebras) yield
large classes of extended Lax pairs $\left(3\right) .$ A general
statement leading to such situations is given by the
Adler-Kostant-Symes theorem.
\begin{Theo}
Let $\mathcal{L}$ be a Lie algebra paired with itself via a
nondegenerate, ad-invariant bilinear form $\langle $ $,$ $\rangle
$, $\mathcal{L}$ having a vector space decomposition
$\mathcal{L}=\mathcal{K}+\mathcal{N}$ with $\mathcal{K}$ and
$\mathcal{N}$ Lie subalgebras. Then, with respect to $\langle $
$,$ $\rangle $, we have the splitting
$\mathcal{L}=\mathcal{L}^{*}=\mathcal{K}^{\perp
}+\mathcal{N}^{\perp }$ and $\mathcal{N}^{*}=\mathcal{K}^{\perp }$
paired with $\mathcal{N}$ via an induced form $\left\langle
\left\langle ,\right\rangle \right\rangle $ inherits the coadjoint
symplectic structure of Kostant and Kirillov; its Poisson bracket
between functions $H_{1}$ and $H_{2}$ on $\mathcal{N}^{*}$ reads
$$\left\{ H_{1},H_{2}\right\} \left( a\right) =\left\langle \left\langle a,
\left[ \nabla _{\mathcal{N}^{*}}H_{1},\nabla
_{\mathcal{N}^{*}}H_{2}\right] \right\rangle \right\rangle \text{
},\text{ }a\in \mathcal{N}^{*}.$$ Let $V\subset \mathcal{N}^*$ be
an invariant manifold under the above co-adjoint action of
$\mathcal{N}$ on $\mathcal{N}^{*}$ and let $\mathcal{A}(V)$ be the
algebra of functions defined on a neighborhood of $V$, invariant
under the coadjoint action of $\mathcal{L}$ (which is distinct
from the $\mathcal{N-N}^*$ action). Then the functions $H$ in
$\mathcal{A}(V)$ lead to commuting Hamiltonian vector fields of
the Lax isospectral form
$$\dot a=\left[ a,pr_{\mathcal{K}}(\nabla H)\right]
\text{ },\text{ }pr_{\mathcal{K}}\text{ }projection\text{
}onto\text{ }\mathcal{K}$$
\end{Theo}
This theorem produces hamiltonian systems having many commuting
integrals ; some precise results are known for interesting classes
of orbits in both the case of finite and infinite dimensional Lie
algebras. Any finite dimensional Lie algebra $\mathcal{L}$ with
bracket $\left[ ,\right] $ and killing form $\left\langle
,\right\rangle $ leads to an infinite dimensional formal Laurent
series extension $\mathcal{L}=\sum_{-\infty
}^{N}A_{i}h^{i}:A_{i}\in \mathcal{L}, \text{ }N\in \mathbb{Z}$
free, with bracket $\left[ \sum A_{i}h^{i},\sum B_{j}h^{j}\right]
=\sum_{i,j}\left[ A_{i},B_{j}\right] h^{i+j},$ and ad-invariant,
symmetric forms $\left\langle \sum A_{i}h^{i},\sum
B_{j}h^{j}\right\rangle _{k}=\sum_{i+j=-k} \left\langle
A_{i},B_{j}\right\rangle,$ depending on $k\in \mathbb{Z}$. The
forms $\left\langle ,\right\rangle _{k}$ are non degenerate if
$\left\langle ,\right\rangle $ is so. Let $\mathcal{L}_{p,q}$
$(p\leq q)$ be the vector space of powers of $h$ between $p$ and
$q$ . A first interesting class of problems is obtained by taking
$\mathcal{L}=\mathcal{G}l(n,\mathbb{R})$ and by putting the form
$\left\langle ,\right\rangle _{1}$ on the Kac-Moody extension.
Then we have the decomposition into Lie subalgebras
$\mathcal{L}=\mathcal{L}_{0,\infty }+\mathcal{L}_{-\infty
,-1}=\mathcal{K}+\mathcal{N}$ with $\mathcal{K=K}^{\perp },$
$\mathcal{N=N}^{\perp }$ and $\mathcal{K=N}^{*}$. Consider the
invariant manifold $V_{m}$ $,$ $m\geq 1$ in $\mathcal{K=N}^{*}$ ,
defined as
$$V_{m}=\left\{ A=\sum_{i=1}^{m-1}A_{i}h^{i}+\alpha h^{m}\text{ },
\text{ }\alpha =diag(\alpha _{1},\cdots ,\alpha _{n})\text{
fixed}\right\} ,$$ with $diag\left( A_{m-1}\right) =0.$
\begin{Theo}
The manifold $V_{m}$ has a natural symplectic structure, the
functions $H=\left\langle f(Ah^{-j}),h^{k}\right\rangle _{1}$ on
$V_{m}$ for good functions $f$ lead to complete integrable
commuting hamiltonian systems of the form
$$\dot{A}=\left[ A,pr_{\mathcal{K}}(f^{\prime }(Ah^{-j})h^{k-j})\right],
\text{ }A=\sum_{i=0}^{m-1}A_{i}h^{i}+\alpha h,$$ and their
trajectories are straight line motions on the jacobian of the
curve $\mathcal{C}$ of genus $\left( n-1\right) \left( nm-2\right)
/2$ defined by $P\left( z,h\right) =\det \left( A-zI\right) =0.$
The coefficients of this polynomial provide the orbit invariants
of $V_{m}$ and an independent set of integrals of the motion (of
particular interest are the flows where $j=m,k=m+1$ which have the
following form
$$\dot{A}=\left[ A\text{ },\text{ }ad_{\beta
\text{ }}ad_{\alpha }^{-1}A_{m-1}+\beta h\right] ,\text{ }\beta
_{i}=f^{\prime } \left( \alpha _{i}\right),$$ the flow depends on
$f$ through the relation $\beta _{i}=f^{\prime }\left( \alpha
_{i}\right) $ only).
\end{Theo}
Another class is obtained by choosing any semi-simple Lie algebra
$L$ . Then the Kac-Moody extension $\mathcal{L}$ equipped with the
form $\left\langle ,\right\rangle =\left\langle ,\right\rangle
_{0}$ has the natural level decomposition $\mathcal{L}=\sum_{i\in
\mathbb{Z}}L_{i},\left[ L_{i,}L_{j}\right] \subset L_{i+j}, \text{
}\left[ L_{0},L_{0}\right] =0,\text{ }L_{i}^{*}=L_{-i}.$ Let
$B^{+}=\sum_{i\geq 0}L_{i}$ and $B^{-}=\sum_{i\langle 0}L_{i}$ .
Then the product Lie algebra $\mathcal{L\times L}$ has the
following bracket and pairing
$$\left[ \left( l_{1},l_{2}\right) ,(l_{1}^{^{\prime }},l_{2}^{^{\prime }})\right]
=\left([l_{1},l_{1}^{^{\prime}}],-[l_{2},l_{2}^{^{\prime}}]\right),
\quad \left\langle \left( l_{1},l_{2}\right) ,(l_{1}^{^{\prime
}},l_{2}^{^{\prime }})\right\rangle =\langle l_{1},l_{1}^{^{\prime
}}\rangle -\langle l_{2},l_{2}^{^{\prime }}\rangle.$$ It admits
the decomposition into $\mathcal{K}+\mathcal{N}$ with
$$\mathcal{K}=\left\{ (l,-l):l\in \mathcal{L}\right\}, \quad
\mathcal{K}^{\perp }=\left\{ (l,l):l\in \mathcal{L}\right\},$$
$$\mathcal{N}=\left\{ (l_{-},l_{+}):l_{-}\in B^{-},l_{+}\in
B^{+},pr_{0}(l_{-})=pr_{0}(l_{+})\right\},$$$$ \mathcal{N}^{\perp
}=\left\{ (l_{-},l_{+}):l_{-}\in B^{-},l_{+}\in
B^{+},pr_{0}(l_{+}+l_{-})=0\right\},$$ where $pr_{0}$ denotes
projection onto $L_{0}.$ Then from the last theorem , the orbits
in $\mathcal{N}^{*}\mathcal=K^{\perp }$ possesses a lot of
commuting hamiltonian vector fields of Lax form:
\begin{Theo}
The N-invariant manifolds $V_{-j,k}=\sum_{-j\leq i\leq
k}L_{i}\subseteq \mathcal{L\simeq K}^{\perp},$ has a natural
symplectic structure and the functions $H(l_{1},l_{2})=f(l_{1})$
on $V_{-j,k}$ lead to commuting vector fields of the Lax form
$$\dot{l}=\left[ l,(pr^{+}-\frac{1}{2}pr_{0})\nabla H\right],
\text{ }pr^{+}\text{ }projection\text{ }onto\text{}B^{+},$$ their
trajectories are straight line motions on the Jacobian of a curve
defined by the characteristic polynomial of elements in
$V_{-j,k}.$
\end{Theo}

Using the van Moerbeke-Mumford linearization method
$\left[21\right] $, Adler and van Moerbeke $\left[1\right] $
showed that the linearized flow could be realized on the jacobian
variety $Jac(\mathcal{C})$ (or some sub-abelian variety of it) of
the algebraic curve (spectral curve) $\mathcal{C}$ associated to
$\left(3\right) $. We then construct an algebraic map from the
complex invariant manifolds of these hamiltonian systems to the
jacobian variety $Jac(\mathcal{C})$ of the curve $\mathcal{C}.$
Therefore all the complex flows generated by the constants of the
motion are straight line motions on these jacobian varieties i.e.
the linearizing equations are given by
$$\int_{s_{1}(0)}^{s_{1}(t)}\omega _{k}+
\int_{s_{2}(0)}^{s_{2}(t)}\omega _{k}+\cdots +
\int_{s_{g}(0)}^{s_{g}(t)}\omega _{k}=c_{k}t\text{ },\text{ }0\leq
k\leq g,$$ where $\omega _{1},\ldots ,\omega _{g}$ span the
$g$-dimensional space of holomorphic differentials on the curve
$\mathcal{C}$ of genus $g.$ In an unifying approach, Griffiths
$\left[8\right] $ has found necessary and sufficient conditions on
$B$ for the Lax flow $\left(3\right) $ to be linearizable on the
jacobi variety of its spectral curve, without reference to
Kac-Moody Lie algebras.\\

Next I schall discuss a number of integrable hamiltonian systems.

\subsection{The Euler rigid body motion.} It express the free
motion of a rigid body around a fixed point. Let $M=\left(
m_{1},m_{2},m_{3}\right) $ be the angular momentum, $\Omega
=\left( m_{1}/I_{1},m_{2}/I_{2},m_{3}/I_{3}\right) $ the angular
velocity and $I_{1},I_{2}$ et $I_{3},$ the principal moments of
inertia about the principal axes of inertia. Then the motion of
the body is governed by
\begin{equation}\label{eqn:euler}
\dot M=M\wedge \Omega.
\end{equation}
If one identifies vectors in $\mathbb{R}^{3}$ with skew-symmetric
matrices by the rule
$$a=\left( a_{1},a_{2},a_{3}\right),\quad
{A}=\left(\begin{array}{ccc}
0&-a_{3}&a_{2}\\
a_{3}&0&-a_{1}\\
-a_{2}&a_{1}&0
\end{array}\right),
$$
then $a\wedge b\longmapsto \left[ A,B\right] =AB-BA.$ Using this
isomorphism between $(\mathbb{R}^{3},\wedge)$ and $(so(3),[,]),$
we write (4) as $\dot M=\left[ M,\Omega \right] ,$ where
$$M= \left(\begin{array}{ccc}
0&-m_{3}&m_{2}\\
m_{3}&0&-m_{1}\\
-m_{2}&m_{1}&0
\end{array}\right)\in so\left( 3\right) ,
\quad \Omega = \left(\begin{array}{ccc}
0&-\omega_{3}&\omega_{2}\\
\omega_{3}&0&-\omega_{1}\\
-\omega_{2}&\omega_{1}&0
\end{array}\right)\in so\left( 3\right) ,
$$
Now $M=I\Omega ,$ this implies that
\begin{equation}\label{eqn:euler}
\dot M=\left[ M,\Lambda M\right] ,
\end{equation}
where
$$\Lambda M=
\left(\begin{array}{ccc}
0&-\lambda _{3}m_{3}&\lambda _{2}m_{2}\\
\lambda _{3}m_{3}&0&-\lambda _{1}m_{1}\\
-\lambda _{2}m_{2}&\lambda _{1}m_{1}&0
\end{array}\right)\in so\left( 3\right),
$$
with $\lambda _{i}\equiv I_{i}^{-1}.$ Equation (5) is explicitly
given by
\begin{eqnarray}
\dot m_{1}&=&\left( \lambda _{3}-\lambda _{2}\right)
m_{2}m_{3},\nonumber\\
\dot m_{2}&=&\left( \lambda _{1}-\lambda _{3}\right) m_{1}m_{3},\\
\dot m_{3}&=&\left( \lambda _{2}-\lambda _{1}\right)
m_{1}m_{2},\nonumber
\end{eqnarray}
and can be written as a hamiltonian vector field
$$\dot x=J\frac{\partial H}{\partial x},\text{ }x=\left( m_{1},m_{2},m_{3}\right) ^{\intercal
},$$ with the hamiltonian $H=\frac{1}{2}\left( \lambda
_{1}m_{1}^{2}+\lambda _{2}m_{2}^{2}+\lambda _{3}m_{3}^{2}\right)
,$ and
$$J=\left(\begin{array}{ccc}
0&-m_{3}&m_{2}\\
m_{3}&0&-m_{1}\\
-m_{2}&m_{1}&0
\end{array}\right)\in so\left(3\right).
$$
We have $\det $ $J=0,$ so $m=2n+k$ and $m-k=rk$ $J.$ Here $m=3$
and $rk$ $J=2$, then $n=k=1.$ The system (6) has beside the energy
$H_{1}=H,$ a trivial invariant $H_{2}$, i.e., such that:
$J\frac{\partial H_{2}}{\partial x}=0,$ or
$$\left(\begin{array}{ccc}
0&-m_{3}&m_{2}\\
m_{3}&0&-m_{1}\\
-m_{2}&m_{1}&0
\end{array}\right)
\left(\begin{array}{c}
\frac{\partial H_{2}}{\partial m_{1}}\\
\frac{\partial H_{2}}{\partial m_{2}}\\
\frac{\partial H_{2}}{\partial m_{3}}
\end{array}\right)
= \left(\begin{array}{c}
0\\
0\\
0
\end{array}\right),
$$
implying $\frac{\partial H_{2}}{\partial m_{1}}=m_{1},\text{
}\frac{\partial H_{2}}{\partial m_{2}}=m_{2},\text{
}\frac{\partial H_{2}}{\partial m_{3}}=m_{3},$ and consequently
$$H_{2}=\frac{1}{2}\left( m_{1}^{2}+m_{2}^{2}+m_{3}^{2}\right) .$$
The system evolves on the intersection of the sphere $H_{1}=c_{1}$
and the ellipsoid $H_{2}=c_{2}.$ In $\mathbb{R}^{3},$ this
intersection will be isomorphic to two circles $\left( \text{with
}\frac{c_{2}}{\lambda _{3}}< c_{1}< \frac{c_{2}}{\lambda
_{1}}\right).$ We shall show that the problem can be integrated in
terms of elliptic functions, as Euler discovered using his then
newly invented theory of elliptic integrals. Observe that the
first equation of (6) reads
\begin{equation}\label{eqn:euler}
\frac{dm_{1}}{m_{2}m_{3}}=\left( \lambda _{3}-\lambda
_{2}\right)dt,
\end{equation}
where $m_1,m_2$ and $m_3$ are related by
 $$\lambda_{1}m_{1}^{2}+\lambda _{2}m_{2}^{2}+\lambda _{3}m_{3}^{2}=c_{1},
 \quad m_{1}^{2}+m_{2}^{2}+m_{3}^{2}=c_{2}.$$ Therefore, if $\lambda
_{2}\neq \lambda _{3},$ we have
$$m_{2}=\pm \sqrt{\frac{c_{2}\lambda _{3}-c_{1}+
\left( \lambda _{1}-\lambda _{3}\right) m_{1}^{2}}{\lambda _{3}-
\lambda _{2}}},\quad m_{3}=\pm \sqrt{\frac{c_{1}-c_{2}\lambda
_{2}+\left( \lambda _{2}-\lambda _{1}\right) m_{1}^{2}}{\lambda
_{3}-\lambda _{2}}}.$$ Substituting these expressions into (7), we
find after integration that the system (6) amounts to an elliptic
integral
$$\int_{m_{1}\left( 0\right) }^{m_{1}\left( t\right) }\frac{dm}
{\sqrt{\left( m^{2}+a\right) \left( m^{2}+b\right)}}=ct,$$ with
respect to the elliptic curve
\begin{equation}\label{eqn:euler}
\mathcal{C}: w^{2}=\left( z^{2}+a\right) \left( z^{2}+b\right),
\end{equation}
with $a=\frac{c_{2}\lambda _{3}-c_{1}}{\lambda _{1}-\lambda
_{3}},\text{ }b =\frac{c_{1}-c_{2}\lambda _{2}}{\lambda
_{2}-\lambda _{1}},\text{ }c =\sqrt{\left( \lambda _{1}-\lambda
_{3}\right) \left( \lambda _{2} -\lambda _{1}\right)}.$ Then the
functions $m_i(t)$ can be expressed in terms of theta-functions
of $t,$ according to the classical inversion of abelian integrals.\\

We shall use the Lax representation of the equations of motion to
show that the linearized Euler flow can be realized on an elliptic
curve isomorphic to the original elliptic curve (8). The solution
to $\left(5\right) $ has the form
$$M\left( t\right) =O\left( t\right) M\left( t\right) M^{\top }
\left( t\right),$$ where $O\left( t\right) $ is one parameter
sub-group of $SO\left(3\right).$ So the hamiltonian flow $\left(
5\right) $ preserves the spectrum of $X$ and therefore its
characteristic polynomial $\det \left( M-zI\right) =-z\left(
z^{2}+m_{1}^{2}+m_{2}^{2}+m_{3}^{2}\right).$ Unfortunately, the
spectrum of a $3\times 3$ skew- symmetric matrix provides only one
piece of information; the conservation of energy does not appear
as part of the spectral information. Therefore one is let to
considering another formulation. The basic observation, due to
Manakov $\left[20\right] ,$ is that equation $\left(5\right) $ is
equivalent to the Lax equation
$$\dot A=\left[ A,B\right] ,$$
where $A=M+\alpha h,\quad B=\Lambda M+\beta h,$ with a formal
indeterminate $h$ and
$${\alpha}=\left(\begin{array}{ccc}
\alpha _{1}&0&0\\
0&\alpha _{2}&0\\
0&0&\alpha _{3}
\end{array}\right),\quad
{\beta}=\left(\begin{array}{ccc}
\alpha _{1}&0&0\\
0&\beta _{2}&0\\
0&0&\beta _{3}
\end{array}\right),
$$
$$\lambda _{1}=\frac{\beta _{3}-\beta _{2}}{\alpha _{3}-\alpha _{2}},
\quad\lambda _{2}=\frac{\beta _{1}-\beta _{3}}{\alpha _{1}-\alpha
_{3}},\quad \lambda _{3}=\frac{\beta _{2}-\beta _{1}}{\alpha
_{2}-\alpha _{1}},$$ and all $\alpha_i$ distinct. The
characteristic polynomial of $A$ is
\begin{eqnarray}
P\left( h,z\right) &=&\det \left( A-zI\right) ,\nonumber\\
&=&\det \left( M+\alpha h-zI\right) ,\nonumber\\
&=&\prod_{j=1}^{3}\left( \alpha_{j}h-z\right) +\left(
\sum_{j=1}^{3}\alpha_{j}m_{j}^{2}\right) h-\left(
\sum_{j=1}^{3}m_{j}^{2}\right) z.\nonumber
\end{eqnarray}
The spectrum of the matrix $A=M+\alpha h$ as a function of $h\in
\mathbb{C}$ is time independent and is given by the zeroes of the
polynomial $P\left( h,z\right),$ thus defining an algebraic curve
(spectral curve). Letting $w=h/z,$ we obtain the following
elliptic curve
$$z^{2}\prod_{j=1}^{3}\left( \alpha _{j}w-1\right)
+2H_{1}w-2H_{2}=0,$$ which is shown to be isomorphic to the
original elliptic curve. Finally, we have the

\begin{Theo}
The Euler rigid body motion is a completely integrable system and
the linearized flow can be realized on an elliptic curve.
\end{Theo}

\subsection{The geodesic flow for a left invariant metric on $SO\left( 4\right).$}
Consider the group $SO(4)$ and its Lie algebra $so(4)$ paired with
itself, via the customary inner product $\left\langle
X,Y\right\rangle =-\frac{1}{2}\text{ }tr\text{ } \left(
X.Y\right),$ where
$$X=
\left(\begin{array}{cccc}
0&-x_{3}&x_{2}&-x_{4}\\
x_{3}&0&-x_{1}&-x_{5}\\
-x_{2}&x_{1}&0&-x_{6}\\
x_{4}&x_{5}&x_{6}&0
\end{array}\right)\in so(4).
$$
A left invariant metric on $SO(4)$ is defined by a non-singular
symmetric linear map $\Lambda :so(4)\longrightarrow so(4),\text{ }
X\longmapsto \Lambda .X,$ and by the following inner product;
given two vectors $gX$ and $gY$ in the tangent space $SO(4)$ at
the point $g\in SO(4),$ $\left\langle gX,gY\right\rangle
=\left\langle X, \Lambda ^{-1}.Y\right\rangle.$ Then the geodesic
flow for this metric takes the following commutator form
(Euler-Arnold equations) :
\begin{equation}\label{eqn:euler}
\dot X=\left[ X,\Lambda .X\right],
\end{equation}
where
$$\Lambda .X=\left(\begin{array}{cccc}
0&-\lambda _{3}x_{3}&\lambda _{2}x_{2}&-\lambda _{4}x_{4}\\
\lambda _{3}x_{3}&0&-\lambda _{1}x_{1}&-\lambda _{5}x_{5}\\
-\lambda _{2}x_{2}&\lambda _{1}x_{1}&0&-\lambda _{6}x_{6}\\
\lambda _{4}x_{4}&\lambda _{5}x_{5}&\lambda _{6}x_{6}&0
\end{array}\right)\in so(4).
$$
In view of the isomorphism between $\left( \mathbb{R}^{6},\wedge
\right) ,$ and $\left( so\left( 4\right) ,\left[ ,\right] \right)
$ we write the system (9) as
\begin{eqnarray}
\dot x_{1}&=&\left( \lambda _{3}-\lambda _{2}\right)
x_{2}x_{3}+\left( \lambda _{6}-\lambda _{5}\right)
x_{5}x_{6},\nonumber\\
\dot x_{2}&=&\left( \lambda _{1}-\lambda _{3}\right)
x_{1}x_{3}+\left( \lambda _{4}-\lambda _{4}\right) x_{4}x_{6},\nonumber\\
\dot x_{3}&=&\left( \lambda _{2}-\lambda _{1}\right)
x_{1}x_{2}+\left( \lambda _{5}-\lambda _{4}\right) x_{4}x_{5},\nonumber\\
\dot x_{4}&=&\left( \lambda _{3}-\lambda _{5}\right)
x_{3}x_{5}+\left( \lambda _{6}-\lambda _{2}\right) x_{2}x_{6},\nonumber\\
\dot x_{5}&=&\left( \lambda _{4}-\lambda _{3}\right)
x_{3}x_{4}+\left( \lambda _{1}-\lambda _{6}\right) x_{1}x_{6},\nonumber\\
\dot x_{6}&=&\left( \lambda _{2}-\lambda _{4}\right)
x_{2}x_{4}+\left( \lambda _{5}-\lambda _{1}\right)
x_{1}x_{5}.\nonumber
\end{eqnarray}
These equations can be written as a hamiltonian vector field
\begin{equation}\label{eqn:euler}
\dot x(t)=J\frac{\partial H}{\partial x},\text{ }x\in
\mathbb{R}^{6},
\end{equation}
with
$$H=\frac{1}{2}\left\langle X,\Lambda X\right\rangle =
\frac{1}{2}\left( \lambda _{1}x_{1}^{2}+\lambda _{2}x_{2}^{2}
+\cdots +\lambda _{6}x_{6}^{2}\right),$$ the hamiltonian and
$$
{J}=\left(\begin{array}{cccccc}
0&-x_{3}&x_{2}&0&-x_{6}&x_{5}\\
x_{3}&0&-x_{1}&x_{6}&0&-x_{4}\\
-x_{2}&x_{1}&0&-x_{5}&x_{4}&0\\
0&-x_{6}&x_{5}&0&-x_{3}&x_{2}\\
x_{6}&0&-x_{4}&x_{3}&0&-x_{1}\\
-x_{5}&x_{4}&0&-x_{2}&x_{1}&0
\end{array}\right)\in so(6).
$$
We have $\det J=0,$ so $m=2n+k$ and $m-k=rk$ $J.$ Here $m=6$ and
$rg$ $J=4$, then $n=k=2.$ The system (10) has beside the energy
$H_1=H,$ two trivial constants of motion :
\begin{eqnarray}
H_{2}&=&\frac{1}{2}\left( x_{1}^{2}+x_{2}^{2}+\cdots
+x_{6}^{2}\right) ,\nonumber\\
H_{3}&=&x_{1}x_{4}+x_{2}x_{5}+x_{3}x_{6}.\nonumber
\end{eqnarray}
Recall that $H_2$ and $H_3$ are called trivial invariants (or
Casimir functions) because $J\frac{\partial H_{2}}{\partial x}=$
$J\frac{\partial H_{3}}{\partial x}=0.$ In order that the
hamiltonian system (10) be completely integrable, it is suffices
to have one more integral, which we take of the form
$$H_{4}=\frac{1}{2}\left( \mu _{1}x_{1}^{2}+\mu _{2}x_{2}^{2}
+\cdots +\mu _{6}x_{6}^{2}\right).$$ The four invariants must be
functionally independent and in involution, so in particular
$$\left\{ H_{4},H_{3}\right\} =\left\langle \frac{\partial H_{4}}
{\partial x},J\frac{\partial H_{3}}{\partial x}\right\rangle=0,$$
i.e.,
\begin{eqnarray}
&&\left( \left( \lambda _{3}-\lambda _{2}\right) \mu _{1}+\left(
\lambda _{1}-\lambda _{3}\right) \mu _{2}+\left( \lambda
_{2}-\lambda _{1}\right) \mu _{3}\right)
x_{1}x_{2}x_{3}\nonumber\\
&&+\left( \left( \lambda _{6}-\lambda _{5}\right) \mu _{1}+ \left(
\lambda _{1}-\lambda _{6}\right) \mu _{5}+ \left( \lambda
_{5}-\lambda _{1}\right)
\mu _{6}\right) x_{1}x_{5}x_{6}\nonumber\\
&&+\left( \left( \lambda _{4}-\lambda _{6}\right) \mu _{2}+ \left(
\lambda _{6}-\lambda _{2}\right) \mu _{4}+ \left( \lambda
_{2}-\lambda _{4}\right)
\mu _{6}\right) x_{2}x_{4}x_{6}\nonumber\\
&&+\left( \left( \lambda _{5}-\lambda _{4}\right) \mu _{3}+\left(
\lambda _{3}-\lambda _{5}\right) \mu _{4}+\left( \lambda
_{4}-\lambda _{3}\right) \mu _{5}\right)
x_{3}x_{4}x_{5}=0.\nonumber
\end{eqnarray}
Then
$$\left( \lambda _{3}-\lambda _{2}\right) \mu _{1}+\left( \lambda _{1}-\lambda _{3}\right) \mu _{2}+\left( \lambda _{2}-\lambda _{1}\right) \mu
_{3}=0,$$
$$\left( \lambda _{6}-\lambda _{5}\right) \mu _{1}+\left( \lambda _{1}-\lambda _{6}\right) \mu _{5}+\left( \lambda _{5}-\lambda _{1}\right) \mu
_{6}=0,$$
$$\left( \lambda _{4}-\lambda _{6}\right) \mu _{2}+\left( \lambda _{6}-\lambda _{2}\right) \mu _{4}+\left( \lambda _{2}-\lambda _{4}\right) \mu
_{6}=0,$$
$$\left( \lambda _{5}-\lambda _{4}\right) \mu _{3}+\left( \lambda _{3}-\lambda _{5}\right) \mu _{4}+\left( \lambda _{4}-\lambda _{3}\right) \mu
_{5}=0.$$ Put
$${\mathcal{A}}=\left(\begin{array}{cccccc}
\lambda _{3}-\lambda _{2}&\lambda _{1}-\lambda _{3}
&\lambda _{2}-\lambda _{1}&0&0&0\nonumber\\
\lambda _{6}-\lambda _{5}&0&0&0&\lambda _{1}-\lambda _{6}
&\lambda _{5}-\lambda _{1}\nonumber\\
0&\lambda _{4}-\lambda _{6}&0&\lambda _{6}-\lambda _{2}
&0&\lambda _{2}-\lambda _{4}\nonumber\\
0&0&\lambda _{5}-\lambda _{4}&\lambda _{3}-\lambda _{5}&\lambda
_{4}-\lambda _{3}&0\nonumber
\end{array}\right).
$$
The number of solutions of this system is equal to the number of
columns of the matrix $\mathcal{A}$ minus the rank of
$\mathcal{A}.$ If $rk\mathcal{A}=4,$ we have two solutions : $\mu
_{i}=1$ lead to the invariant $H_{2}$ and $\mu _{i}=\lambda _{i}$
lead to the invariant $H_{3}.$ This is unacceptable. If
$rk\mathcal{A}=3,$ each four-order minor of $\mathcal{A}$ is
singular. Now
$$\left(\begin{array}{cccc}
\lambda _{3}-\lambda _{2}&\lambda _{1}-\lambda _{3}
&\lambda _{2}-\lambda _{1}&0\nonumber\\
\lambda _{6}-\lambda _{5}&0&0&0\nonumber\\
0&\lambda _{4}-\lambda _{6}&0&\lambda _{6}-\lambda _{2}\nonumber\\
0&0&\lambda _{5}-\lambda _{4}&\lambda _{3}-\lambda _{5}\nonumber
\end{array}\right)=-\left( \lambda _{6}-\lambda _{5}\right) C,
$$
$$\left(\begin{array}{cccc}
\lambda _{1}-\lambda _{3}&\lambda _{2}-\lambda _{1}
&0&0\nonumber\\
0&0&0&\lambda _{1}-\lambda _{6}\nonumber\\
\lambda _{4}-\lambda _{6}&0&\lambda _{6}-\lambda _{2}&0\nonumber\\
0&\lambda _{5}-\lambda _{4}&\lambda _{3}-\lambda _{5}&\lambda
_{4}-\lambda _{3}\nonumber
\end{array}\right)=\left( \lambda _{1}-\lambda _{6}\right) C,
$$
$$\left(\begin{array}{cccc}
\lambda _{2}-\lambda _{1}&0&0&0\nonumber\\
0&0&\lambda _{1}-\lambda _{6}&\lambda _{5}-\lambda _{1}\nonumber\\
0&\lambda _{6}-\lambda _{2}&0&\lambda _{2}-\lambda _{4}\nonumber\\
\lambda _{5}-\lambda _{4}&\lambda _{3}-\lambda _{5}&\lambda
_{4}-\lambda _{3}&0\nonumber
\end{array}\right)=-\left( \lambda _{2}-\lambda _{1}\right) C,
$$
where
\begin{eqnarray}
C\equiv &&\lambda _{1}\lambda _{6}\lambda _{4}+\lambda _{1}\lambda
_{2}\lambda _{5}-\lambda _{1}\lambda _{2}\lambda _{4}+\lambda
_{3}\lambda _{6}\lambda _{5}-\lambda _{3}\lambda _{6}\lambda
_{4}-\lambda _{3}\lambda _{2}\lambda _{5}\nonumber\\
&&+\lambda_{4}\lambda _{2}\lambda _{5}+\lambda _{4}\lambda
_{1}\lambda _{3}-\lambda _{4}\lambda _{1}\lambda _{5}+\lambda
_{6}\lambda _{2}\lambda _{3}-\lambda _{6}\lambda _{2}\lambda
_{5}-\lambda _{1}\lambda _{6}\lambda _{3},\nonumber
\end{eqnarray}
and it follows that the condition for which these minors are zero
is $C=0.$ Notice that this relation holds by cycling the indices :
$1\curvearrowright 4, 2\curvearrowright 5,
\curvearrowright3\curvearrowright 6.$ Under Manakov $[20]$
conditions,
\begin{eqnarray}
\lambda _{1}&=&\frac{\beta _{2}-\beta _{3}}{\alpha _{2}-\alpha
_{3}},\text{ }\lambda _{2}=\frac{\beta _{1}-\beta _{3}}{\alpha
_{1}-\alpha _{3}},\text{ }\lambda _{3}=\frac{\beta _{1}-\beta
_{2}}{\alpha _{1}-\alpha _{2}},\\
\lambda _{4}&=&\frac{\beta _{1}-\beta _{4}}{\alpha _{1}-\alpha
_{4}},\text{ }\lambda _{5}=\frac{\beta _{2}-\beta _{4}}{\alpha
_{2}-\alpha _{4}},\text{ }\lambda _{6}=\frac{\beta _{3}-\beta
_{4}}{\alpha _{3}-\alpha _{4}},\nonumber
\end{eqnarray}
where $\alpha _{i},\beta _{i}$ $\in \mathbb{C},$
$\prod_{i<j}\left( \alpha _{i}-\beta _{j}\right) \neq 0,$
equations (10) admits a Lax equation with an indeterminate $h:$
\begin{eqnarray}
&&\dot {( \overbrace{X+\alpha h})}=\left[ X+\alpha h,\Lambda
X+\beta h\right] ,\\
{\alpha}&=&\left(\begin{array}{cccc}
\alpha _{1}&0&0&0\nonumber\\
0&\alpha _{2}&0&0\nonumber\\
0&0&\alpha _{3}&0\nonumber\\
0&0&0&\alpha _{4}\nonumber
\end{array}\right),\quad
{\beta}=\left(\begin{array}{cccc}
\beta_{1}&0&0&0\nonumber\\
0&\beta _{2}&0&0\nonumber\\
0&0&\beta _{3}&0\nonumber\\
0&0&0&\beta _{4}\nonumber
\end{array}\right).
\end{eqnarray}
$$\Updownarrow$$
\begin{eqnarray}
\dot X&=&\left[ X,\Lambda .X\right] \Leftrightarrow \left(
9\right),\nonumber\\
\left[ X,\beta \right]&+&\left[ \alpha,\Lambda .X\right] =0
\Leftrightarrow \left(11\right),\nonumber\\
\left[ \alpha ,\beta \right]&=&0\text{ $\text{trivially satisfied
for diagonal matrices.}$}\nonumber
\end{eqnarray}
The parameters $\mu _{1},\ldots ,\mu _{6}$ can be parameterized
(like $\lambda _{1},\ldots ,\lambda _{6}$) by :
$$\mu _{1}=\frac{\gamma _{2}-\gamma _{3}}{\alpha _{2}-\alpha _{3}},\text{ }\mu _{2}=
\frac{\gamma _{1}-\gamma _{3}}{\alpha _{1}-\alpha _{3}},\text{
}\mu _{3}= \frac{\gamma _{1}-\gamma _{2}}{\alpha _{1}-\alpha
_{2}}$$
$$\mu _{4}=\frac{\gamma _{1}-\gamma _{4}}{\alpha _{1}-\alpha _{4}},\text{ }\mu _{5}=
\frac{\gamma _{2}-\gamma _{4}}{\alpha _{2}-\alpha _{4}},\text{
}\mu _{6}= \frac{\gamma _{3}-\gamma _{4}}{\alpha _{3}-\alpha
_{4}}.$$ To use the method of isospectral deformations, consider
the Kac-Moody extension $\left( n=4\right)$: $\mathcal{L}=\left\{
\sum_{-\infty }^{N}A_{i}h^{i}: N\text{arbitrary }\in
\mathbb{Z},\text{ } A_{i}\in gl(n,\mathbb{R})\right\},$ of $gl(n,
\mathbb{R})$ with the bracket: $\left[ \sum A_{i}h^{i},\sum
B_{j}h^{j}\right] =\sum_{k}\left( \sum_{i+j=k}\left[
A_{i},B_{j}\right] \right)h^{k},$ and the ad-invariant form:
$\left\langle \sum A_{i}h^{i},\sum B_{j}h^{j}\right\rangle
=\sum_{i+j=-1}\langle A_{i}, B_{j}\rangle,$ where
$\langle,\rangle$ is the usual form defined on $gl(n,
\mathbb{R}).$ Let $\mathcal{K}$ and $\mathcal{N}$ be respectively
the $\geq 0$ and $<0$ powers of $h$ in $\mathcal{L},$ then
$\mathcal{L}=\mathcal{K}+\mathcal{N},$ for the pairing defined
above $\mathcal{K}=\mathcal{K}^{\perp },\text{ }\mathcal{N}
=\mathcal{N}^{\perp},$ so that $\mathcal{K}=\mathcal{N}^*$. The
orbits described in this way come equipped with a symplectic
structure with Poisson bracket $\left\{ H_{1},H_{2}\right\} \left(
\alpha \right) =\left\langle \alpha ,\left[ \nabla
_{\mathcal{K}^{*}}H_{1},\nabla _{\mathcal{K}^{*}}H_{2}\right]
\right\rangle ,$ where $\alpha \in \mathcal{K}^{*}$\ and $\nabla
_{\mathcal{K}^{*}}H\in \mathcal{K}.$ According to the
Adler-Kostant-Symes theorem, the flow $\left(12\right) $\ is
hamiltonian on an orbit through the point $X+ah,$ $X\in so(4)$)
formed by the coadjoint action of the subgroup
$G_{\mathcal{N}}\subset SL\left( n\right) $ of lower triangular
matrices on the dual Kac-Moody algebra $\mathcal{N}^{*}\approx
\mathcal{K}^{\perp }=\mathcal{K}$. As a consequence, the
coefficients of $z^{i}h^{i}$\ appearing in curve :
\begin{equation}\label{eqn:euler}
\Gamma :\text{ }\left\{ \left( z,h\right) \in \mathbb{C}^{2}:\det
\left( X+ah-zI\right) =0\right\},
\end{equation}
associated to the equation $\left(12\right) ,$\ are invariant of
the system in involution for the symplectic structure of this
orbit. Notice that
$$
\det \left( gXg^{-1}\right) =\det X =\left(
x_{1}x_{4}+x_{2}x_{5}+x_{3}x_{6}\right) ^{2},$$
$$ tr\left(gXg^{-1}\right) ^{2}=tr\left( gX^{2}g^{-1}\right) =tr\left(
X^{2}\right)= -2\left( x_{1}^{2}+x_{2}^{2}+\cdots
+x_{6}^{2}\right).$$ Also the complex flows generated by these
invariants can be realized as straight lines on the abelian
variety defined by the periods of curve $\Gamma.$ Explicitly,
equation $\left(13\right) $ looks as follows
$$
\Gamma :\quad  \prod_{i=1}^{4}\left( \alpha _{i}h-z\right) +2H_{4}
h^{2}-2H_{1}zh+2H_{2} z^{2}+H_{3}^{2}=0,
$$
where $H_{1}(X)=c_{1},\text{ }H_{2}(X)=c_{2},\text{
}H_{3}(X)=2H=c_{3},\text{ }H_{4}(X)=c_{4}.$ with
$c_{1},c_{2},c_{3},c_{4}$ generic constants. $\Gamma $\ is a curve
of genus $3$\ and it has a natural involution $\sigma :\Gamma
\rightarrow \Gamma \text{ },\text{ }(z,h)\mapsto (-z,-h).$
Therefore the jacobian variety $Jac(\Gamma )$ of $\Gamma $ splits
up into an even and old part : the even part is an elliptic curve
$\Gamma _{0}$\ = $\Gamma /\sigma $\ and the odd part is a
$2-$dimensional abelian surface $Prym(\Gamma /\Gamma _{0})$ called
the Prym variety : $Jac(\Gamma )=\Gamma _{0}+Prym(\Gamma /\Gamma_0
)$. The van Moerbeke-Mumford linearization method provides then an
algebraic map from the complex affine variety
$\bigcap_{i=1}^{4}\left\{ H_{i}(X)=c_{i}\right\} \subset
\mathbb{C}^{6}$ to the Jacobi variety $Jac(\Gamma ).$\ By the
antisymmetry of $\Gamma $\ $,$\ this map sends this variety to the
Prym variety $Prym(\Gamma /\Gamma_0 ):$
$$\bigcap_{i=1}^{4}\left\{ H_{i}(X)=c_{i}\right\} \rightarrow
Prym(\Gamma /\Gamma_0 )\text{ },\text{ }p\mapsto
\sum_{k=1}^{3}s_{k},$$ and the complex flows generated by the
constants of the motion are straight lines on $Prym(\Gamma
/\Gamma_0 )$. Finally, we have the

\begin{Theo}
The geodesic flow $\left(9\right) $ is a hamiltonian system with
$$H\equiv H_{1}=\frac{1}{2}\left( \lambda _{1}x_{1}^{2}+
\lambda _{2}x_{2}^{2}+\cdots +\lambda _{6}x_{6}^{2}\right),$$ the
hamiltonian. It has two trivial invariants
\begin{eqnarray}
H_{2}&=&\frac{1}{2}\left( x_{1}^{2}+x_{2}^{2}+\cdots
+x_{6}^{2}\right) ,\nonumber\\
H_{3}&=&x_{1}x_{4}+x_{2}x_{5}+x_{3}x_{6}.\nonumber
\end{eqnarray}
Moreover, if
\begin{eqnarray}
&&\lambda _{1}\lambda _{6}\lambda _{4}+\lambda _{1}\lambda
_{2}\lambda _{5}-\lambda _{1}\lambda _{2}\lambda _{4}+\lambda
_{3}\lambda _{6}\lambda _{5}-\lambda _{3}\lambda _{6}\lambda
_{4}-\lambda _{3}\lambda _{2}\lambda _{5}\nonumber\\
&&+\lambda_{4}\lambda _{2}\lambda _{5}+\lambda _{4}\lambda
_{1}\lambda _{3}-\lambda _{4}\lambda _{1}\lambda _{5}+\lambda
_{6}\lambda _{2}\lambda _{3}-\lambda _{6}\lambda _{2}\lambda
_{5}-\lambda _{1}\lambda _{6}\lambda _{3}=0,\nonumber
\end{eqnarray}
the system (9) has a fourth independent constant of the motion of
the form
$$H_{4}=\frac{1}{2}\left( \mu _{1}x_{1}^{2}+\mu _{2}x_{2}^{2}+
\cdots +\mu _{6}x_{6}^{2}\right).$$ Then the system (9) is
completely integrable and can be linearized on the Prym variety
$Prym(\Gamma /\Gamma_0 )$.
\end{Theo}

\subsection{The Toda lattice.}
The Toda lattice equations (discretized version of the Korteweg-de
Vries equation\footnote{In short K-dV equation : $\frac{\partial
u} {\partial t}-6u\frac{\partial u}{\partial x}+\frac{\partial
^{3}u}{\partial x^{3}}=0.$ This is an infinite-dimensional
completely integrable system.}) motion of $n$ particles with
exponential restoring forces are governed by the following
hamiltonian
$$H=\frac{1}{2}\sum_{i=1}^{N}p_{i}^{2}+\sum_{i=1}^{N}e^{q_{i}-q_{i+1}},
\qquad q_{N+1}=q_1.$$ The hamiltonian equations can be written as
follows
$$\dot q_{i}=p_{i},\quad \dot p_{i}=-e^{q_{i}-q_{i+1}}
+e^{q_{i-1}-q_{i}}.$$ In term of the Flaschka's variables
$\left[6\right] :$ $a_{i}=\frac{1}{2}e^{q_{i}-q_{i+1}},$
$b_{i}=-\frac{1}{2}p_{i},$ Toda's equations take the following
form
\begin{equation}\label{eqn:euler}
\dot a_{i}=a_{i}\left( b_{i+1}-b_{i}\right),\quad \dot
b_{i}=2(a_{i}^2-a_{i-1}^2),
\end{equation}
with $b_{N+1}=b_{1}$ and $a_{0}=a_{N}.$ To show that the system
(14) is completely integrable, one should find $N$ first integrals
independent and in involution each other. From the second
equation, we have
$$\frac{d}{dt}\sum_{i=1}^{N}b_{i}=\sum_{i=1}^{N}\frac{db_{i}}{dt}=0,$$
and we normalize the $b_i$'s by requiring that
$\sum_{i=1}^{N}b_{i}=0.$ This is a first integral for the system.
We further define $N\times N$ matrices $A$ and $B$ with
$$
{A}=\left(\begin{array}{ccccc}
b_{1}&a_{1}&0&\cdots&a_{N}\\
a_{1}&b_{2}&\vdots &&\vdots \\
0&\ddots&\ddots&\ddots&0\\
\vdots&&\ddots&b_{N-1}&a_{N-1}\\
a_{N}&\cdots&0&a_{N-1}&b_{N}
\end{array}\right), {B}=\left(\begin{array}{ccccc}
0&a_{1}&\cdots&\cdots&-a_{N}\\
-a_{1}&0&\vdots &&\vdots \\
\vdots&\ddots&\ddots&\ddots&\vdots\\
\vdots&&\ddots&\ddots&a_{N-1}\\
a_{N}&\cdots&\cdots&-a_{N-1}&0
\end{array}\right).
$$
Then (14) is equivalent to the Lax equation
$$\dot A=\left[ B,A\right].$$
From theorem 3, we know that the quantities
$I_k=\frac{1}{k}trA^{k},$ $1\leq k\leq N,$ are first integrals of
motion : To be more precise
$$\dot I_k=tr(\dot A.A^{k-1})=tr([B,A].A^{k-1})
=tr(BA^k-ABA^{k-1})=0.$$ Notice that $I_1$ is the first integral
already know. Since these $N$ first integrals are shown to be
independent and in involution each other, the system (14) is thus
completely integrable.

\subsection{The Garnier potential.}
Consider the hamiltonian
\begin{equation}\label{eqn:euler}
H=\frac{1}{2}\left( x_{1}^{2}+x_{2}^{2}\right) -\frac{1}{2}\left(
\lambda _{1}y_{1}^{2}+\lambda _{2}y_{2}^{2}\right)
+\frac{1}{4}\left( y_{1}^{2}+y_{2}^{2}\right) ^{2},
\end{equation}
where $\lambda _{1}$ and $\lambda _{2}$ are constants. The
corresponding system is given by
\begin{eqnarray}
\dot y_{1}&=&x_{1},\qquad\dot x_{1}=\left( \lambda
_{1}-y_{1}^{2}-y_{2}^{2}\right) y_{1},\\
\dot y_{2}&=&x_{2},\qquad\dot x_{2}=\left( \lambda
_{2}-y_{2}^{2}-y_{1}^{2}\right) y_{2}.\nonumber
\end{eqnarray}

\begin{Theo}
The system $\left( 16\right) $ has the additional first integral
\begin{eqnarray}
H_{2}&=&\frac{1}{4}\left( \left( x_{1}y_{2}-x_{2}y_{1}\right)
^{2}-\left( \lambda _{2}y_{1}^{4}+\lambda _{1}y_{2}^{4}\right)
-\left( \lambda _{1}+\lambda _{2}\right) y_{1}^{2}y_{2}^{2}\right)
\nonumber\\
&&+\frac{1}{2}\left( \lambda _{1}\lambda _{2}\left(
y_{1}^{2}+y_{2}^{2}\right) -\left( \lambda _{2}x_{1}^{2}+\lambda
_{1}x_{2}^{2}\right) \right) .\nonumber
\end{eqnarray}
and is completely integrable. The flows generated by $H_{1}=H$(15)
and $H_{2}$ are straight line motions on the jacobian variety of a
smooth genus two hyperelliptic curve $\mathcal{H}$(17) associated
to a Lax equation.
\end{Theo}
\emph{Proof}. We consider the Lax representation in the form $\dot
A=\left[ A,B\right],$ with the following ansatz for the Lax
operator
$${A}=\left(\begin{array}{cc}
U&V\\
W&-U
\end{array}\right),\qquad
{B}=\left(\begin{array}{cc}
0&1\\
R&0
\end{array}\right)$$
where
\begin{eqnarray}
V&=&-(h-\lambda_{1})(h-\lambda_{2})(1+\frac{1}{2}(\frac{y_{1}^{2}}
{h-\lambda_{1}}+\frac{y_{2}^{2}}{h-\lambda_{2}})) ,\nonumber\\
U&=&\frac{1}{2}(h-\lambda_{1})(h-\lambda_{2})(
\frac{x_{1}y_{1}}{h-\lambda_{1}}+\frac{x_{2}y_{2}}{h-\lambda_{2}}) ,\nonumber\\
W&=&(h-\lambda_{1})(h-\lambda_{2})(\frac{1}{2}(
\frac{x_{1}^{2}}{h-\lambda_{1}}+\frac{x_{2}^{2}}{h-\lambda_{2}})
-h+\frac{1}{2}(y_{1}^{2}+y_{2}^{2})),\nonumber\\
R&=&h-y_{1}^{2}-y_{2}^{2}.\nonumber
\end{eqnarray}
We form the curve in $\left( z,h\right) $ space
$$P(h,z)=\det \left( A-zI\right) =0,$$
whose coefficients are functions of the phase space. Explicitly,
this equation looks as follows
\begin{eqnarray}
\mathcal{H}: z^{2}&=&P_{5}\left( h\right),\\
&=&(h-\lambda_{1})(h-\lambda_{2})(h^{3}-(\lambda_{1}+\lambda_{2})h^{2}+(
\lambda_{1}\lambda_{2}-H_{1})h-H_{2}),\nonumber
\end{eqnarray}
with $H_1$ (15) the hamiltonian and a second quartic integral
$H_{2}$ of the form
\begin{eqnarray}
H_{2}&=&-\frac{1}{4}(\lambda_{2}y_{1}^{4}+\lambda_{1}y_{2}^{4}+(\lambda_{1}
+\lambda_{2})y_{1}^{2}y_{2}^{2}-(x_{1}y_{2}-x_{2}y_{1})^{2})\nonumber\\
&&-\frac{1}{2}(\lambda_{2}x_{1}^{2}+\lambda_{1}x_{2}^{2}-\lambda_{1}\lambda_{2}(
y_{1}^{2}+y_{2}^{2})).\nonumber
\end{eqnarray}
The functions $H_1$ and $H_2$ commute : $\left\{
H_{1},H_{2}\right\}=0$ and the system (16) is completely
integrable. The curve $\mathcal{H}$ determined by the fifth-order
equation $\left(17\right) $ is smooth, hyperelliptic and its genus
is $2.$ Obviously, $\mathcal{H}$ is invariant under the
hyperelliptic involution $\left( h,z\right) \curvearrowright
\left( h,-z\right) .$ Using the van Moerbeke-Mumford linearization
method, we show that the linearized flow could be realized on the
jacobian variety $Jac\left( \mathcal{H} \right) $ of the genus 2
curve $\mathcal{H}.$ For generic $c=\left( c_{1},c_{2}\right) \in
\mathbb{C}^{2}$ the affine variety defined by
$$M_{c}=\overset{2}{\underset{i=1}{\bigcap }}\left\{ x\in \mathbb{C}^{4}:H_{i}\left( x\right)
=c_{i}\right\},$$ is a smooth affine surface. According to the
schema of $\left[5\right],$ we introduce coordinates $s_{1}$ and
$s_{2}$ on the surface $M_{c}$, such that $M_{c}\left(
s_{i}\right) =0,\lambda _{1}\neq \lambda _{2},$ i.e.,
$$s_{1}+s_{2}=\frac{1}{2}\left( y_{1}^{2}+y_{2}^{2}\right) +\lambda _{1}+
\lambda _{2},\quad s_{1}s_{2}=\frac{1}{2}\left( \lambda
_{2}y_{1}^{2}+\lambda _{1}y_{2}^{2}\right) +\lambda _{1}\lambda
_{2}.$$ After some algebraic manipulations, we obtain the
following equations for $s_{1}$ and $s_{2}:$
$$\dot s_{1}=2\frac{\sqrt{P_{5}\left( s_{1}\right)
}}{s_{1}-s_{2}},\quad \dot s_{2}=2\frac{\sqrt{P_{5}\left(
s_{2}\right) }}{s_{2}-s_{1}},$$ where $P_{5}\left( s\right) $ is
defined by $\left(17\right).$ These equations can be integrated by
the abelian mapping
$$\mathcal{H}\longrightarrow Jac\left( \mathcal{H}\right) =
\mathbb{C}^{2}/L\text{ },\text{ }p\longmapsto \left(
\int_{p_{0}}^{p}\omega _{1}\text{ }, \text{
}\int_{p_{0}}^{p}\omega _{2}\right) ,$$ where the hyperelliptic
curve $\mathcal{H}$ of genus two is given by the equation $(17),$
$L$ is the lattice generated by the vectors $n_{1}+\Omega
n_{2},(n_{1},n_{2}) \in \mathbb{Z}^{2},\Omega$\ is the matrix of
period of the curve $\mathcal{H}$, $(\omega _{1},\omega _{2})$\ is
a canonical basis of holomorphic differentials on $\mathcal{H},$\
i.e.,
$$\omega _{1}=\frac{ds}{\sqrt{P_{5}\left( s\right) }}\text{ },\quad
\omega _{2}=\frac{sds}{\sqrt{P_{5}\left( s\right) }},$$ and
$p_{0}$\ is a fixed point. This concludes the proof of the
theorem.

\subsection{The coupled nonlinear Schr\"{o}dinger equations.}
The system of two coupled nonlinear Schr\"{o}dinger equations is
given by
\begin{eqnarray}
i\frac{\partial a}{\partial z}+\frac{\partial ^{2}a}{\partial
t^{2}}+\Omega _{0}a +\frac{2}{3}\left( \left| a\right| ^{2}+\left|
b\right| ^{2}\right) a+
\frac{1}{3}\left( a^{2}+b^{2}\right) \overline{a}&=&0,\\
i\frac{\partial b}{\partial z}+\frac{\partial ^{2}b}{\partial
t^{2}}-\Omega _{0}b +\frac{2}{3}\left( \left| a\right| ^{2}+\left|
b\right| ^{2}\right) b+ \frac{1}{3}\left( a^{2}+b^{2}\right)
\overline{b}&=&0,\nonumber
\end{eqnarray}
where $a\left( z,t\right) $ and $b\left( z,t\right) $ are
functions of $z$ and $t,$ the bar ``$-$'' denotes the complex
conjugation, ``$\left| {}\right| $'' denotes the modulus and
$\Omega _{0}$ is a constant. These equations play a significant
role in mathematics, with a important number of physical
applications. We seek solutions of $\left(18\right)$ in the
following form
$$a\left( z,t\right) =y_{1}\left( t\right) \exp \left( i\Omega z\right) ,\quad
b\left( z,t\right) =y_{2}\left( t\right) \exp \left( i\Omega
z\right) ,$$ where $y_{1}\left( t\right) $ et $y_{2}\left(
t\right) $ are two functions and $\Omega $\ is an arbitrary
constant. Then we obtain the system
$$\ddot y_{1}+\left( y_{1}^{2}+y_{2}^{2}\right) y_{1}=
\left( \Omega -\Omega _{0}\right) y_{1},$$
$$\ddot y_{2}+\left( y_{1}^{2}+y_{2}^{2}\right) y_{2}=
\left( \Omega +\Omega _{0}\right) y_{2}.$$ The latter coincides
obviously with $\left(16\right) $ for $\lambda _{1}=\Omega -\Omega
_{0}$ and $\lambda _{2}=\Omega +\Omega _{0}.$

\subsection{The Yang-Mills equations.}
We consider the Yang-Mills system for a field with gauge group
$SU(2):$
$$\triangledown _{j}F_{jk}=\frac{\partial F_{jk}}{\partial \tau _{j}}+
\left[ A_{j},F_{jk}\right] =0,$$ where $F_{jk}, A_{j}\in
T_{e}SU(2), 1\leq j,k\leq 4$ and $F_{jk}=\frac{\partial
A_{k}}{\partial \tau _{j}}- \frac{\partial A_{j}}{\partial \tau
_{k}}+\left[ A_{j},A_{k}\right].$ The self-dual Yang-Mills (SDYM)
equations is an universal system for which some reductions include
all classical tops from Euler to Kowalewski (0+1-dimensions),
K-dV, Nonlinear Schrödinger, Sine-Gordon, Toda lattice and N-waves
equations (1+1-dimensions), KP and D-S equations (2+1-dimensions).
In the case of homogeneous double-component field, we have
$\partial_{j}A_{k}=0,j\neq 1 , A_{1}=A_{2}=0, A_{3}=n_{1}U_{1}\in
su( 2) , A_{4}=n_{2}U_{2}\in su(2) $ where $n_{i}$ are $su(2)
$-generators (i.e., they satisfy commutation relations :
$n_1=[n_2,[n_1,n_2]], n_2=[n_1,[n_2,n_1]]$). The system becomes
$$\frac{\partial ^{2}U_{1}}{\partial t^{2}}+U_{1}U_{2}^{2}=0,\qquad
\frac{\partial ^{2}U_{2}}{\partial t^{2}}+U_{2}U_{1}^{2}=0,$$ with
$t=\tau _{1}.$ By setting $U_{j}=q_{j},$ $\frac{\partial
U_{j}}{\partial t}=p_{j},$ $j=1,2,$ Yang-Mills equations are
reduced to hamiltonian system
$$\dot x=J\frac{\partial H}{\partial x},\text{ \quad }x=
(q_{1},q_{2},p_{1},p_{2})^{\intercal },\text{ \quad }
{J}=\left(\begin{array}{cc}
O&-I\\
I&O
\end{array}\right),
$$ with $H=\frac{1}{2}( p_{1}^{2}+p_{2}^{2}+q_{1}^{2}q_{2}^{2}),$ the hamiltonian. The
symplectic transformation $p_{1}\curvearrowleft \frac{\sqrt{2}}{2}
( p_{1}+p_{2}) , p_{2}\curvearrowleft \frac{\sqrt{2}}{2} (
p_{1}-p_{2}), q_{1}\curvearrowleft \frac{1}{2}( \root{4}\of{2})(
q_{1}+iq_{2}) , q_{2}\curvearrowleft \frac{1}{2}( \root{4}\of{2})
( q_{1}-iq_{2}),$ takes this hamiltonian into
$$
H=\frac{1}{2}(p_{1}^{2}+p_{2}^{2})+\frac{1}{4}q_{1}^{4}
+\frac{1}{4}q_{2}^{4}+\frac{1}{2}q_{1}^{2}q_{2}^{2},$$ which
coincides with $\left(15\right) $ for $\lambda _{1}=\lambda
_{2}=0.$
\section{Algebraic complete integrability}
We give some results about abelian surfaces which will be used, as
well as the basic techniques to study two-dimensional algebraic
completely integrable systems. Let $M=\mathbb{C}/\Lambda$ be a
$n-$dimensional abelian variety where $\Lambda$ is the lattice
generated by the $2n$ columns $\lambda_1,\ldots,\lambda_{2n}$ of
the $n\times 2n$ period matrix $\Omega$ and let $D$ be a divisor
on $M.$ Define $\mathcal{L}(\mathcal{D})=\{f \mbox{meromorphic on}
M  : (f)\geq-\mathcal{D}\},$ i.e., for $\mathcal{D}=\sum
k_j\mathcal{D}_j$ a function $f\in\mathcal{L}(\mathcal{D})$ has at
worst a $k_j-$fold pole along $\mathcal{D}_j.$ The divisor
$\mathcal{D}$ is called ample when a basis $(f_0,\ldots,f_N)$ of
$\mathcal{L}(k\mathcal{D})$ embeds $M$ smoothly into
$\mathbb{P}^N$ for some $k,$ via the map $M\rightarrow
\mathbb{P}^N,\text{ }p\mapsto [1:f_{1}(p):...:f_{N}(p)],$ then
$k\mathcal{D}$ is called very ample. It is known that every
positive divisor $\mathcal{D}$ on an irreducible abelian variety
is ample and thus some multiple of $\mathcal{D}$ embeds $M$ into
$\mathbb{P}^N.$ By a theorem of Lefschetz, any $k\geq 3$ will
work. Moreover, there exists a complex basis of $\mathbb{C}^n$
such that the lattice expressed in that basis is generated by the
columns of the $n\times 2n$ period matrix
$$ \left(\begin{array}{ccccc}
\delta_1&&0&|&\\
&\ddots&&|&Z\\
0&&\delta_n&|&
\end {array}\right),$$
with $Z^\top =Z, \mbox{Im}Z>0, \delta_j\in \mathbb{N}^*$ and
$\delta_j|\delta_{j+1}.$ The integers $\delta_j$ which provide the
so-called polarization of the abelian variety $M$ are then related
to the divisor as follows : \begin{equation}\label{eqn:euler} \dim
\mathcal{L}(\mathcal{D})=\delta_1\ldots \delta_n.
\end{equation}
In the case of a $2-$dimensional abelian varieties (surfaces),
even more can be stated : the geometric genus $g$ of a positive
divisor $\mathcal{D}$ (containing possibly one or several curves)
on a surface $M$ is given by the adjunction formula
\begin{equation}\label{eqn:euler}
g(\mathcal{D})=\frac{K_{M}.\mathcal{D}+\mathcal{D}.\mathcal{D}}{2}+1,
\end{equation}
where $K_M$ is the canonical divisor on $M,$ i.e., the zero-locus
of a holomorphic $2-$form, $\mathcal{D}.\mathcal{D}$ denote the
number of intersection points of $\mathcal{D}$ with
$a+\mathcal{D}$ (where $a+\mathcal{D}$ is a small translation by
$a$ of $\mathcal{D}$ on $M$), where as the Riemann-Roch theorem
for line bundles on a surface tells you that
\begin{equation}\label{eqn:euler}
\chi(\mathcal{D})=p_{a}(M)+1+\frac{1}{2}(\mathcal{D}.\mathcal{D}-\mathcal{D}K_M),
\end{equation}
where $p_{a}\left( M\right) $ is the arithmetic genus of $M$ and
$\chi(\mathcal{D})$ the Euler characteristic of $\mathcal{D}.$ To
study abelian surfaces using Riemann surfaces on these surfaces,
we recall that
\begin{eqnarray}
\chi( \mathcal{D}) &=&\dim {H}^{0}(
M,\mathcal{O}_{M}(\mathcal{D}))-\dim {H}^{1}(M,\mathcal{O}_{M}(D)),\nonumber\\
&=&\dim{\mathcal{L}}(\mathcal{D})-\dim {H}^{1}(M,\Omega ^{2}(
\mathcal{D}\otimes K_M^*)),\mbox{(Kodaira-Serre duality)},\nonumber\\
&=&\dim {\mathcal{L}}(\mathcal{D}), \mbox{(Kodaira vanishing
theorem)},
\end{eqnarray}
whenever $\mathcal{D}\otimes K_M^*$ defines a positive line
bundle. However for abelian surfaces, $K_M$ is trivial and
$p_a(M)=-1;$ therefore combining relations (19), (20), (21) and
(22),
$$
\chi \left( \mathcal{D}\right)=\dim
{\mathcal{L}}(\mathcal{D})=\frac{\mathcal{D}.\mathcal{D}}{2} =
g\left( \mathcal{D}\right) -1=\delta_1 \delta_2.
$$
A divisor $\mathcal{D}$ is called projectively normal, when the
natural map $\mathcal{L}(\mathcal{D})^{\otimes k}\rightarrow
\mathcal{L}(k\mathcal{D}),$ is surjective, i.e., every function of
$\mathcal{L}(k\mathcal{D})$ can be written as a linear combination
of k-fold products of functions of $\mathcal{L}(\mathcal{D}).$ Not
every very ample divisor $\mathcal{D}$ is projectively normal but
if $\mathcal{D}$ is linearly equivalent to $k\mathcal{D}_0$ for
$k\geq 3$ for some divisor $\mathcal{D}_0,$ then $\mathcal{D}$ is
projectively normal.\\
Now consider the exact sheaf sequence
$$0\longrightarrow \mathcal{O}_C \overset{\pi^* }{\longrightarrow
}\mathcal{O}_{\widetilde{C}}\longrightarrow X\longrightarrow 0,$$
where $C$ is a singular connected Riemann surface,
$\widetilde{C}=\sum C_j$ the corresponding set of smooth Riemann
surfaces after desingularization and $\pi :
\widetilde{C}\rightarrow C$ the projection. The exactness of the
sheaf sequence shows that the Euler characteristic
$$
\mathcal{X}(\mathcal{O})=\dim H^0(\mathcal{O})-\dim
H^1(\mathcal{O}),
$$
satisfy
\begin{equation}\label{eqn:euler}
\mathcal{X}(\mathcal{O}_C)-\mathcal{X}(\mathcal{O}_{\widetilde{C}})+\mathcal{X}(X)=0,
\end{equation}
where $\mathcal{X}(X)$ only accounts for the singular points p of
$C;$ $\mathcal{X}(X_p)$ is the dimension of the set of holomorphic
functions on the different branches around p taken separately,
modulo the holomorphic functions on the Riemann surface $C$ near
that singular point. Consider the case of a planar singularity (in
this paper, we will be concerned by a tacnode for which
$\mathcal{X}(X)=2,$ as well), i.e., the tangents to the branches
lie in a plane. If $f_{j}(x,y)=0$ denote the $j^{th}$ branch of
$C$ running through p with local parameter $s_j,$ then
$$\mathcal{X}(X_p)
=\dim{\Pi_{j}\mathbb{C}[[s_j]]}/\frac{\mathbb{C}[[x,y]]}{\Pi_{j}f_{j}(x,y)}.$$
So using (22) and Serre duality, we obtain
$\mathcal{X}(\mathcal{O}_C)=1-g(C)$ and
$\mathcal{X}(\mathcal{O}_{\widetilde{C}})=n-\sum_{j=1}^{n}g(C_j).$
Also, replacing in the formula (23), gives
$$
g(C)=\sum_{j=1}^{n}g(C_j)+\mathcal{X}(X)+1-n.
$$
Finally, recall that a K\"{a}hler variety is a variety with a
K\"{a}hler metric, i.e., a hermitian metric whose associated
differential $2$-form of type $(1,1)$ is closed. The complex torus
$\mathbb{C}^{2}/lattice$ with the euclidean metric $\sum
dz_{i}\otimes d\overline{z}_{i}$ is a K\"{a}hler variety and any
compact complex variety that can be embedded in projective space
is also a K\"{a}hler variety. Now, a compact complex K\"{a}hler
variety having as many independent meromorphic functions as its
dimension is a projective variety.

Consider now hamiltonian problems of the form
\begin{equation}\label{eqn:euler}
X_{H}:\dot {x}=J\frac{\partial H}{\partial x}\text{ }\equiv
f(x),\text{ }x\in \mathbb{R}^{m},
\end{equation}
where $H$\ is the hamiltonian and $J=J(x)$\ is a skew-symmetric
matrix with polynomial entries in $x,$\ for which the
corresponding Poisson bracket $\{H_{i},H_{j}\}=\langle
\frac{\partial H_i}{\partial x},J\frac{\partial H_j} {\partial
x}\rangle,$ satisfies the Jacobi identities. The system (24) with
polynomial right hand side will be called algebraic complete
integrable (a.c.i.) when :\\
$i)$\ The system possesses $n+k$\ independent polynomial
invariants $H_{1},...,H_{n+k}$ of which $k$ lead to zero vector
fields $J\frac{\partial H_{n+i}}{\partial x}\left( x\right) =0,$
$1\leq i\leq k,$ the $n$ remaining ones are in involution (i.e.,
$\left\{ H_{i},H_{j}\right\} =0$) and $m=2n+k.$\ For most values
of $c_{i}\in \mathbb{R},$ the invariant varieties
$\overset{n+k}{\underset{i=1}{\bigcap }}\left\{ x\in
\mathbb{R}^{m}:H_{i}=c_{i}\right\} $\ are assumed compact and
connected. Then, according to the Arnold-Liouville theorem, there
exists a diffeomorphism
$$\overset{n+k}{\underset{i=1}{\bigcap }}\left\{ x\in
\mathbb{R}^{m}:H_{i}=c_{i}\right\} \rightarrow
\mathbb{R}^{n}/Lattice,$$and the solutions of the system (24) are
straight lines motions on these tori.\\
$ii)$\ The invariant varieties, thought of as affine varieties in
$\mathbb{C}^{m}$ can be completed into complex algebraic tori,
i.e.,$$ \overset{n+k}{\underset{i=1}{\bigcap }}\{H_{i}=c_{i},x\in
\mathbb{C}^{m}\}\cup \mathcal{D}=\mathbb{C}^{n}/Lattice, $$ where
$\mathbb{C}^{n}/Lattice$ is a complex algebraic torus (i.e.,
abelian variety) and $\mathcal{D}$ a divisor. Algebraic means that
the torus\ can be defined as an intersection
$\displaystyle{\bigcap_{i=1}^{M}\{ P_{i}(X_{0},...,X_{N})=0\}} $
involving a large number of homogeneous polynomials $P_{i}.$\ In
the natural coordinates $(t_{1},...,t_{n})$\ of
$\mathbb{C}^{n}$/$Lattice$\ coming from $\mathbb{C}^{n},$\ the
functions $x_{i}=x_{i}(t_{1},...,t_{n})$\ are meromorphic and (24)
defines straight line motion on $\mathbb{C}^{n}/Lattice.$\
Condition $i)$\ means, in particular, there is an algebraic map
$(x_{1}(t),...,x_{m}(t) ) \mapsto (\mu_{1}(t) ,...,\mu_{n}(t))$
making the following sums linear in $t$\
:$$\sum_{i=1}^{n}\int_{\mu_{i}(0)}^{\mu_{i}(t)}\omega _{j}=d
_{j}t\text{ },\text{ }1\leq j\leq n,\text{ }d _{j}\in
\mathbb{C},$$where $\omega _{1},...,\omega _{n}$\ denote
holomorphic differentials on some algebraic curves.\\
The existence of a coherent set of Laurent solutions :
\begin{equation}\label{eqn:euler}
x_{i}=\sum_{j=0}^{\infty }x_{i}^{(j) }t^{j-k_{i}},\text{\quad
}k_{i}\in \mathbb{Z},\text{ \quad some }k_{i}> 0,
\end{equation}
depending on $dim\ (phase$\ $space)-1=$\ $m-1$\ free parameters is
necessary and sufficient for a hamiltonian system with the right
number of constants of motion to be a.c.i. So, if the hamiltonian
flow (24) is a.c.i., it means that the variables $x_{i}$\ are
meromorphic on the torus $\mathbb{C}^{n}/Lattice$\ and by
compactness they must blow up along a codimension one subvariety
(a divisor) $\mathcal{D}\subset \mathbb{C}^{n}/Lattice.$\ By the
a.c.i. definition, the flow (24) is a straight line motion in
$\mathbb{C}^{n}/Lattice$\ and thus it must hit the divisor
$\mathcal{D}$\ in at least one place. Moreover through every point
of $\mathcal{D},$ there is a straight line motion and therefore a
Laurent expansion around that point of intersection. Hence the
differential equations must admit Laurent expansions which depend
on the $n-1$\ parameters defining $\mathcal{D}$\ and the $n+k$\
constants $c_{i}$\ defining the torus $\mathbb{C}^{n}/Lattice$\ ,
the total count is therefore
$m-1=dim\ (phase\ space)-1$\ parameters. \\
Assume now hamiltonian flows to be (weight)-homogeneous with a
weight $\nu _{i}\in \mathbb{N},$ going with each variable $x_{i},$
i.e.,
$$f_{i}\left( \alpha ^{\nu _{1}}x_{1},...,\alpha ^{\nu _{m}}x_{m}\right)
=\alpha ^{\nu _{i}+1}f_{i}\left( x_{1},...,x_{m}\right) , \text{
}\forall \alpha \in \mathbb{C}.$$ Observe that then the constants
of the motion $H$ can be chosen to be (weight)-homogeneous :
$$H\left( \alpha ^{\nu _{1}}x_{1},...,\alpha ^{\nu _{m}}x_{m}\right)=
\alpha ^{k}H\left( x_{1},...,x_{m}\right) ,\text{ }k\in
\mathbb{Z}.$$ If the flow is algebraically completely integrable,
the differential equations $\left(24\right) $ must admits Laurent
series solutions $\left( 25\right) $ depending on $m-1$ free
parameters. We must have $k_{i}=\nu _{i}$ and coefficients in the
series must satisfy at the 0$^{th}$step non-linear equations,
\begin{equation}\label{eqn:euler}
f_{i}\left( x_{1}^{\left( 0\right) },...,x_{m}^{\left( 0\right)
}\right) +g_{i}x_{i}^{\left( 0\right) }=0,\text{ }1\leq i\leq m,
\end{equation}
and at the k$^{th}$step, linear systems of equations :
\begin{equation}\label{eqn:euler}
\left( L-kI\right) z^{\left( k\right) }=
\left\{\begin{array}{rl} 0&\mbox{ for } k=1\\
\mbox{some polynomial in}& x^{\left( 1\right) },...,x^{\left(
k-1\right)} \mbox{ for } k>1,
\end{array}\right.
\end{equation}
where
$$L=\text{ Jacobian map of }\left(26\right) =
\text{ }\frac{\partial f}{\partial z}+gI\mid _{z=z^{\left(
0\right)}}.$$ If $m-1$ free parameters are to appear in the
Laurent series, they must either come from the non-linear
equations $\left(26\right) $ or from the eigenvalue problem
$\left(27\right) ,$ i.e., $L$ must have at least $m-1$ integer
eigenvalues. These are much less conditions than expected, because
of the fact that the homogeneity $k$ of the constant $H$ must be
an eigenvalue of $L$ Moreover the formal series solutions are
convergent as a consequence of the majorant method. Next we assume
that the divisor is very ample and in addition projectively
normal. Consider a point $p\in\mathcal{D},$ a chart $U_j$ around
$p$ on the torus and a function $y_j$ in
$\mathcal{L}(\mathcal{D})$ having a pole of maximal order at $p.$
Then the vector $(1/y_j, y_1/y_j,\ldots,y_N/y_j)$ provides a good
system of coordinates in $U_j.$ Then taking the derivative with
regard to one of the flows
$$(\frac{y_i}{y_j})\dot{}=\frac{\dot{y_i}y_j-y_i\dot{y_j}}{y_j^2},\quad
1\leq j\leq N,$$ are finite on $U_j$ as well. Therefore, since
$y_j^2$ has a double pole along $\mathcal{D},$ the numerator must
also have a double pole (at worst), i.e.,
$\dot{y_i}y_j-y_i\dot{y_j}\in \mathcal{L}(2\mathcal{D}).$ Hence,
when $\mathcal{D}$ is projectively normal, we have that
$$(\frac{y_i}{y_j})\dot{}=\sum_{k,l}a_{k,l}(\frac{y_k}{y_j})(\frac{y_l}{y_j}),$$
i.e., the ratios $y_i/y_j$ form a closed system of coordinates
under differentiation. At the bad points, the concept of
projective normality play an important role: this enables one to
show that $y_i/y_j$ is a bona fide Taylor series starting from
every point in a neighbourhood of the point
in question.\\

To prove the algebraic complete integrability of a given
hamiltonian system, the main steps of the method are :

- The first step is to show the existence of the Laurent
solutions, which requires an argument precisely every time $k$ is
an integer eigenvalue of $L$ and therefore $L-kI$ is not
invertible.

- One shows the existence of the remaining constants of the motion
in involution so as to reach the number $n+k.$

- For given $c_{1},...,c_{m},$ the set
$$\mathcal{D}\equiv \left\{\begin{array}{rl}
&x_{i}\left( t\right)=t^{-\nu _{i}}\left( x_{i}^{\left( 0\right)
}+x_{i}^{\left( 1\right) }t+x_{i}^{\left(
2\right) }t^{2}+\cdots \right) ,1\leq i\leq m,\\
&\mbox{ Laurent solutions such that }: H_{j}\left( x_{i}\left(
t\right) \right) =c_{j}+\mbox{ Taylor part },
\end{array}\right\}
$$
defines one or several $n-1$ dimentional algebraic varieties
(divisor) having the property that
\begin{eqnarray}
\overset{n+k}{\underset{i=1}{\bigcap }}\left\{ H_{i}=c_{i},z\in
\mathbb{C}^{m}\right\} \cup \mathcal{D}\text{ }&=& \text{a smooth
compact,
connected variety } \nonumber\\
&&\text{with $n$ commuting vector fields }\nonumber\\
&&\text{independent at every point.}\nonumber\\
&=&\text{ a complex algebraic torus
}T^{n}=\mathbb{C}^{n}/Lattice.\nonumber
\end{eqnarray}
The flows $\,J\frac{\partial H_{k+i}}{\partial z}$
$,...,J\frac{\partial H_{k+n}}{\partial z}$ are straight line
motions on $T^{n}.$ \\
From the divisor $\mathcal{D}$, a lot of information can be
obtained with regard to the periods and the action-angle
variables.

\subsection{A five-dimensional system}
Consider the following system of five differential
equations in the unknowns $z_1,\ldots,z_5:$
\begin{eqnarray}
\dot{z}_1&=&2z_4,\quad\qquad \dot{z}_3=z_2(3z_1+8z_2^2) ,\nonumber\\
\dot{z}_2&=&z_3,\quad\qquad
\dot{z}_4=z_{1}^{2}+4z_{1}z_{2}^{2}+z_{5},
\end{eqnarray}
$$\dot{z}_5=2z_{1}z_{4}+4z_2^2z_4-2z_{1}z_{2}z_{3}.$$
The following  three quartics are constants of motion for this
system
\begin{eqnarray}
F_1&=&\frac{1}{2}z_{5}-z_{1}z_{2}^{2}+\frac{1}{2}z_{3}^{2}
-\frac{1}{4}z_{1}^{2}-2z_{2}^{4},\nonumber\\
F_2&=&z_5^2-z_1^2z_5+4z_1z_2z_3z_4-z_1^2z_3^2+\frac{1}{4}z_1^4-4z_2^2z_4^2 ,\\
F_3&=&z_{1}z_{5}+z_{1}^{2}z_{2}^{2}-z_{4}^{2}.\nonumber
\end{eqnarray}
This system is completely integrable and the hamiltonian structure
is defined by the Poisson bracket $\left\{ F,H\right\}
=\left\langle \frac{\partial F}{\partial z}, J\frac{\partial
H}{\partial z}\right\rangle =\sum_{k,l=1}^{5}J_{kl}\frac{\partial
F}{\partial z_{k}}\frac{\partial H}{\partial z_{l}},$ where
$\frac{\partial H}{\partial z}=(\frac{\partial H}{\partial
z_{1}},\frac{\partial H}{\partial z_{2}},\frac{\partial
H}{\partial z_{3}},\frac{\partial H}{\partial
z_{4}},\frac{\partial H}{\partial z_{5}})^\top,$ and
$$J=\left[\begin{array}{ccccc}
0&0&0&2z_1&4z_4\\
0&0&1&0&0\\
0&-1&0&0&-4z_1z_2\\
-2z_1&0&0&0&2z_5-8z_1z_{2}^{2}\\
-4z_4&0&4z_1z_2&-2z_5+8z_1z_{2}^{2}&0
\end{array}\right],$$
is a skew-symmetric matrix for which the corresponding Poisson
bracket satisfies the Jacobi identities. The system (28) can be
written as $\dot{z}=J\frac{\partial H}{\partial z},
z=(z_{1},z_{2},z_{3},z_{4},z_{5})^\top,$ where $H=F_{1}.$ The
second flow commuting with the first is regulated by the equations
$\dot{z}=J\frac{\partial F_{2}}{\partial z},
z=(z_{1},z_{2},z_{3},z_{4},z_{5})^\top.$ These vector fields are
in involution : $\{F_1,F_2\}=\langle \frac{\partial
F_{1}}{\partial z},J\frac{\partial F_{2}}{\partial z}\rangle=0,$
and the remaining one is casimir : $J\frac{\partial
F_{3}}{\partial z}=0.$ The invariant variety $A$ defined by
\begin{equation}\label{eqn:euler}
A=\bigcap_{k=1}^{2}\{z:F_k(z)=c_k\}\subset\mathbb{C}^5,
\end{equation}
is a smooth affine surface for generic values of
$(c_{1},c_{2},c_{3}) \in \mathbb{C}^{3}$. So, the question I
address is how does one find the compactification of $A$ into an
abelian surface? The idea of the direct proof we shall give here
is closely related to the geometric spirit of the (real)
Arnold-Liouville theorem. Namely, a compact complex
$n$-dimensional variety on which there exist $n$ holomorphic
commuting vector fields which are independent at every point is
analytically isomorphic to a $n$-dimensional complex torus
$\mathbb{C}^{n}/Lattice$ and the complex flows generated by the
vector fields are straight lines on this complex torus. Now, the
main problem will be to complete $A\left(30\right)$ into a non
singular compact complex algebraic variety $\widetilde{A}=A\cup
\mathcal{D}$ in such a way that the vector fields $X_{F_1}$ and
$X_{F_2}$ generated respectively by $F_1$ and $F_2,$ extend
holomorphically along a divisor $\mathcal{D}$ and remain
independent there. If this is possible, $\widetilde{A}$ is an
algebraic complex torus (an abelian variety) and the coordinates
$z_{1},\ldots,z_{5}$ restricted to $A$ are abelian functions. A
naive guess would be to take the natural compactification
$\overline{A}$ of $A$ by projectivizing the equations:
$\overline{A}=\bigcap_{k=1}^{3}\{F_{k}(Z)=c_kZ_0^4\} \subset
\mathbb{P}^{5}.$ Indeed, this can never work for a general reason
: an abelian variety $\widetilde{A}$ of dimension bigger or equal
than two is never a complete intersection, that is it can never be
described in some projective space $ \mathbb{P}^{n}$ by $n$-dim
$\widetilde{A}$ global polynomial homogeneous equations. In other
words, if $A$ is to be the affine part of an abelian surface,
$\overline{A}$ must have a singularity somewhere along the locus
at infinity $\overline{A}\cap \left\{ Z_{0}=0\right\} .$ In fact,
we shall show that the existence of meromorphic solutions to the
differential equations (28) depending on 4 free parameters can be
used to manufacture the tori, without ever going through the
delicate procedure of blowing up and down. Information about the
tori can then be gathered from the divisor.

\begin{Theo}
The system (28) possesses Laurent series solutions which depend on
4 free parameters : $\alpha, \beta, \gamma $ and $\theta.$ These
meromorphic solutions restricted to the surface $A$(30) are
parameterized by two copies $\mathcal{C}_{-1}$ and
$\mathcal{C}_{1}$ of the same Riemann surface (32) of genus 7.
\end{Theo}
\emph{Proof}. The first fact to observe is that if the system is
to have Laurent solutions depending on 4 free parameters
$\alpha,\beta,\gamma,\theta$, the Laurent decomposition of such
asymptotic solutions must have the following form
\begin{eqnarray}
z_{1}&=&\frac{1}{t}\alpha -\frac{1}{2}\alpha ^{2}+\beta t
-\frac{1}{16}\alpha \left( \alpha ^{3}+4\beta \right) t^{2}+\gamma t^{3}+\cdots, \nonumber\\
z_{2}&=&\frac{1}{2t}\varepsilon -\frac{1}{4}\varepsilon \alpha
+\frac{1}{8}\varepsilon \alpha ^{2}t
-\frac{1}{32}\varepsilon \left(-\alpha ^{3}+12\beta \right) t^{2}+\theta t^{3}+\cdots, \nonumber\\
z_{3}&=&-\frac{1}{2t^{2}}\varepsilon +\frac{1}{8}\varepsilon
\alpha ^{2}\allowbreak
-\frac{1}{16}\varepsilon \left( -\alpha ^{3}+12\beta \right) t+3\theta t^{2}+\cdots,\\
z_{4}&=&-\frac{1}{2t^{2}}\alpha +\frac{1}{2}\beta
-\frac{1}{16}\alpha \left( \alpha ^{3}+4\beta \right) t+\frac{3}{2}\gamma t^{2}+\cdots, \nonumber\\
z_{5}&=&\frac{1}{2t^{2}}\alpha ^{2}-\frac{1}{4t}\left( \alpha
^{3}+4\beta \right) +\allowbreak \frac{1}{4}\alpha \left( \alpha
^{3}+2\beta \right) -\left( \alpha ^{2}\beta -2\gamma
+4\varepsilon \theta \alpha \right)t +\cdots, \nonumber
\end{eqnarray}
with $\varepsilon=\pm 1.$ Using the majorant method, we can show
that these series are convergent. Substituting the Laurent
solutions $(31)$ into $(29)$: $F_{1}=c_{1},$ $F_{2}=c_{2}$ and
$F_{3}=c_{3},$ and equating the $t^0$-terms yields
\begin{eqnarray}
F_{1}&=&\allowbreak \frac{7}{64}\alpha ^{4}-\frac{1}{8}\alpha
\beta
-\frac{5}{2}\varepsilon \theta=c_1, \nonumber\\
F_{2}&=&\frac{1}{16}\left( 4\beta -\alpha ^{3}\right) \left(
4\alpha ^{2}\beta
-\alpha ^{5}+64\varepsilon \theta \alpha -32\gamma \right)=c_2,\nonumber\\
F_{3}&=&-\frac{1}{32}\alpha ^{6}-\beta ^{2}-\frac{1}{4}\alpha
^{3}\beta -3\varepsilon \theta \alpha ^{2}+4\alpha
\gamma=c_3.\nonumber
\end{eqnarray}
Eliminating $\gamma $ and $\theta$ from these equations, leads to
an equation connecting the two remaining parameters $\alpha$ and
$\beta$ :
\begin{equation}\label{eqn:euler}
\mathcal{C}: 64\beta ^{3}-16\alpha ^{3}\beta ^{2}-4\left( \alpha
^{6}-32\alpha ^{2}c_{1}-16c_{3}\right) \beta $$ $$+\alpha \left(
32c_{2}-32\alpha ^{4}c_{1}+\alpha ^{8}-16\alpha
^{2}c_{3}\right)=0.
\end{equation}
The Laurent solutions restricted to the surface $A$(30) are thus
parameterized by two copies $\mathcal{C}_{-1}$ and
$\mathcal{C}_{1}$ of the same Riemann surface $\mathcal{C}$(32).
According to the Riemann-Hurwitz formula, the genus of the Riemann
surface
$\mathcal{C}$ is 7, which establishes the theorem.\\

In order to embed $\mathcal{C}$ into some projective space, one of
the key underlying principles used is the Kodaira embedding
theorem, which states that a smooth complex manifold can be
smoothly embedded into projective space $\mathbb{P}^N$ with the
set of functions having a pole of order k along positive divisor
on the manifold, provided k is large enough; fortunately, for
abelian varieties, k need not be larger than three according to
Lefshetz. These functions are easily constructed from the Laurent
solutions (31) by looking for polynomials in the phase variables
which in the expansions have at most a k-fold pole. The nature of
the expansions and some algebraic proprieties of abelian varieties
provide a recipe for when to terminate our search for such
functions, thus making the procedure implementable. Precisely, we
wish to find a set of polynomial functions $\{f_0,\ldots,f_N\},$
of increasing degree in the original variables $z_1,\ldots,z_5$
having the property that the embedding $\mathcal{D}$ of
$\mathcal{C}_1+\mathcal{C}_{-1}$ into $\mathbb{P}^N$ via those
functions satisfies the relation : geometric genus
$(\mathcal{D})\equiv g(\mathcal{D})=N+2.$ A this point, it may be
not so clear why $\mathcal{D}$ must really live on an abelian
surface. Let us say, for the moment, that the equations of the
divisor $\mathcal{D}$ (i.e., the place where the solutions blow
up), as a Riemann surface traced on the abelian surface
$\widetilde{A}$ (to be constructed in theorem 12), must be
understood as relations connecting the free parameters as they
appear firstly in the expansions (31). This means that (32) must
be understood as relations connecting $\alpha$ and $\beta.$ Let
$$L^{(r)}=\left\{\begin{array}{rl}
&\mbox{polynomials}\quad f=f(z_,\ldots,z_5)\\& \mbox{of degree}
\leq r, \quad\mbox{such that}\\&
f(z(t))=t^{-1}(z^{(0)}+\ldots),\\
&\mbox{with} \quad z^{(0)}\neq 0 \quad\mbox{on}\quad
\mathcal{D}\\& \mbox{and with}\quad z(t) \quad\mbox{as in}\quad
(4)
\end{array}\right\}/[F_k=c_k, k=1,2,3],$$ and let
$(f_0,f_1,\ldots,f_{N_r})$ be a basis of $L^{(r)}.$ We look for r
such that : $g(\mathcal{D}^{(r)})=N_r+2, \mathcal{D}^{(r)}\subset
\mathbb{P}^{N_r}.$ We shall show that it is unnecessary to go
beyond r=4.

\begin{Theo}
a) The spaces $L^{(r)}$, nested according to weighted degree, are
generated as follows
\begin{eqnarray}
L^{(1)}&=&\{f_0,f_1,f_2\},\nonumber\\
L^{(2)}&=&L^{(1)}\oplus\{f_3,f_4,f_5,f_6\},\nonumber\\
L^{(3)}&=&L^{(2)}\oplus\{f_7,f_8,f_9,f_{10}\},\nonumber\\
L^{(4)}&=&L^{(3)}\oplus\{f_{12},f_{13},f_{14},f_{15}\},
\end{eqnarray}
where $f_0=1, f_1=z_{1}, f_2=z_{2},f_3=2z_{5}-z_{1}^{2},
f_4=z_{3}+2\varepsilon z_{2}^{2}, f_5=z_{4}+\varepsilon
z_{1}z_{2},  f_6=\left[ f_{1},f_{2}\right],
f_7=f_{1}(f_{1}+2\varepsilon f_{4}), f_8=f_{2}(f_{1}+2\varepsilon
f_{4}), f_{9}=z_{4}(f_{3}+2\varepsilon f_{6}),
f_{10}=z_{5}(f_{3}+2\varepsilon f_{6}),
f_{11}=f_{5}(f_{1}+2\varepsilon f_{4}),
f_{12}=f_{1}f_{2}(f_{3}+2\varepsilon f_{6}),
f_{13}=f_{4}f_{5}+\left[ f_{1},f_{4}\right], f_{14}=\left[
f_{1},f_{3}\right]+2\varepsilon \left[ f_{1},f_{6}\right],
f_{15}=f_{3}-2z_{5}+4f_{4}^{2},$ with $[s_j,s_k]=\dot s_j
s_k-s_j\dot s_k,$ the
wronskien of $s_k$ and $s_j$.\\
b) $L^{(4)}$ provides an embedding of $\mathcal{D}^{(4)}$ into
projective space $\mathbb{P}^{15}$ and $\mathcal{D}^{(4)}$ has
genus 17.
\end{Theo}
\emph{Proof}. a) The proof of a) is straightforward and can be
done by inspection of the expansions (31).\\
b) It turns out that neither $L^{(1)},$ nor $L^{(2)},$ nor
$L^{(3)},$ yields a Riemann surface of the right genus; in fact $
g(\mathcal{D}^{(r)})\neq \dim L^{(r)}+1, r=1,2,3.$ For instance,
the embedding into $\mathbb{P}^{2}$ via $L^{(1)}$ does not
separate the sheets, so we proceed to $L^{(2)}$ and the
corresponding embedding into $\mathbb{P}^{6}$ is unacceptable
since $g(\mathcal{D}^{(2)})-2>6$ and $\mathcal{D}^{(2)}\subset
\mathbb{P}^{6}\neq \mathbb{P}^{g-2},$ which contradicts the fact
that $N_r=g(\mathcal{D}^{(2)})-2.$ So we proceed to $L^{(3)}$ and
we consider the corresponding embedding into $\mathbb{P}^{10}$,
according to the functions $(f_0,\ldots,f_{10}).$  For finite
values of $\alpha$ and $\beta,$ dividing the vector
$(f_0,\ldots,f_{10})$ by $f_2$ and taking the limit $t\rightarrow
0,$ to yield $[0:2\varepsilon \alpha:1:-\varepsilon
(4\beta-\alpha^3):-\alpha: -\varepsilon
\alpha^2:\frac{1}{2}(4\beta-\alpha^3):\varepsilon
\alpha^3:\frac{1}{2}\alpha ^2:\frac{1}{4}\varepsilon \alpha
^3(4\beta-\alpha^3):-\frac{1}{4}\varepsilon \alpha
^4(4\beta-\alpha^3)].$ The point $\alpha=0$ require special
attention. Indeed near $\alpha=0,$ the parameter $\beta$ behaves
as follows : $\beta\sim 0,i\sqrt{c_3},-i\sqrt{c_3}.$ Thus near
$(\alpha,\beta)=(0,0),$ the corresponding point is mapped into the
point $[0:0:1:0:0:0:0:0:0:0:0]$ in $\mathbb{P}^{10}$ which is
independent of $\varepsilon=\pm1,$ whereas near the point
$(\alpha,\beta)=(0,i\sqrt{c_3})$ (resp.
$(\alpha,\beta)=(0,-i\sqrt{c_3})$) leads to two different points :
$[0:0:1:-4\varepsilon i\sqrt{c_3}:0:0:2\varepsilon
i\sqrt{c_3}:0:0:0:0]$ (resp. $[0:0:1:4\varepsilon
i\sqrt{c_3}:0:0:-2\varepsilon i\sqrt{c_3}:0:0:0:0]$), according to
the sign of $\varepsilon.$ The Riemann surface (31) has three
points covering $\alpha=\infty,$ at which $\beta$ behaves as
follows : $\beta \sim -\frac{1279}{216}\alpha ^{3},
\frac{1}{432}\alpha ^{3}\left( 1333-1295i\sqrt{3}\right),
\frac{1}{432}\alpha ^{3}\left( 1333+1295i\sqrt{3}\right).$ Then by
dividing the vector $(f_0,\ldots,f_{10})$ by $f_{10},$ the
corresponding point is mapped into the point
$[0:0:0:0:0:0:0:0:0:0:1]$ in $\mathbb{P}^{10}$. Thus,
$g(\mathcal{D}^{(3)})-2>10$ and $\mathcal{D}^{(2)}\subset
\mathbb{P}^{10}\neq \mathbb{P}^{g-2},$ which contradicts the fact
that $N_r=g(\mathcal{D}^{(3)})-2.$ Consider now the embedding
$\mathcal{D}^{(4)}$ into $\mathbb{P}^{15}$ using the 16 functions
$f_0,\ldots,f_{15}$ of $L^{(4)}$(33). It is easily seen that these
functions separate all points of the Riemann surface (except
perhaps for the points at $\alpha=\infty$ and $\alpha=\beta=0$) :
The Riemann surfaces $\mathcal{C}_1$ and $\mathcal{C}_{-1}$ are
disjoint for finite values of $\alpha$ and $\beta$ except for
$\alpha=\beta=0$; dividing the vector $(f_0,\ldots,f_{15})$ by
$f_2$ and taking the limit $t\rightarrow 0,$ to yield
$[0:2\varepsilon \alpha:1:-\varepsilon (4\beta-\alpha^3):-\alpha:
-\varepsilon \alpha^2:\frac{1}{2}(4\beta-\alpha^3):\varepsilon
\alpha^3:\frac{1}{2}\alpha ^2:\frac{1}{4}\varepsilon \alpha
^3(4\beta-\alpha^3):-\frac{1}{4}\varepsilon \alpha
^4(4\beta-\alpha^3):-\frac{1}{2}\varepsilon \alpha ^4:-\frac{1}{4}
\alpha ^3(4\beta-\alpha^3):\frac{3}{4}\alpha \left( 4\beta -\alpha
^{3}\right):\varepsilon \alpha ^{3}\left( 4\beta -\alpha
^{3}\right):-2\varepsilon \alpha ^{3}],$ As before, the point
$\alpha=0$ require special attention and the parameter $\beta$
behaves as follows : $\beta\sim 0,i\sqrt{c_3},-i\sqrt{c_3}.$ Thus
near $(\alpha,\beta)=(0,0),$ the corresponding point is mapped
into the point $[0:0:1:0:0:0:0:0:0:0:0:0:0:0:0:0]$ in
$\mathbb{P}^{15}$ which is independent of $\varepsilon=\pm1,$
whereas near the point $(\alpha,\beta)=(0,i\sqrt{c_3})$ (resp.
$(\alpha,\beta)=(0,-i\sqrt{c_3})$) leads to two different points :
$[0:0:1:-4\varepsilon i\sqrt{c_3}:0:0:2\varepsilon
i\sqrt{c_3}:0:0:0:0:0:0:0:0:0]$ (resp. $[0:0:1:4\varepsilon
i\sqrt{c_3}:0:0:-2\varepsilon i\sqrt{c_3}:0:0:0:0:0:0:0:0:0]$),
according to the sign of $\varepsilon.$ About the point
$\alpha=\infty,$ it is appropriate to divide by $f_{10};$ then the
corresponding point is mapped into the point
$[0:0:0:0:0:0:0:0:0:0:1:0:0:0:0:0],$ in $\mathbb{P}^{15}$ which is
independent of $\varepsilon$. The divisor $\mathcal{D}^{(4)}$
obtained in this way has genus 17 and $\mathcal{D}^{(4)}\subset
\mathbb{P}^{15}=\mathbb{P}^{g-2},$
as desired. This ends the proof of the theorem.\\

Let $\mathcal{L}=L^{(4)}$ and $\mathcal{D}=\mathcal{D}^{(4)}$.
Next we wish to construct a surface strip around $\mathcal{D}$
which will support the commuting vector fields. In fact,
$\mathcal{D}$ has a good chance to be very ample divisor on an
abelian surface, still to be constructed.

\begin{Theo}
The variety $A$(30) generically is the affine part of an abelian
surface $\widetilde{A}$. The reduced divisor at infinity
$\widetilde{A}\backslash A=\mathcal{C}_1+\mathcal{C}_{-1},$
consists of two copies $\mathcal{C}_1$ and $\mathcal{C}_{-1}$ of
the same genus 7 Riemann surface $\mathcal{C}$(32). The system of
differential equations (28) is algebraically completely integrable
and the corresponding flows evolve on $\widetilde{A}.$
\end{Theo}
\emph{Proof}. We need to attaches the affine part of the
intersection of the three invariants $F_1, F_2, F_3$ so as to
obtain a smooth compact connected surface in $\mathbb{P}^{15}.$ To
be precise, the orbits of the vector field (28) running through
$\mathcal{D}$ form a smooth surface $\Sigma$ near $\mathcal{D}$\
such that $\Sigma \backslash A \subseteq \widetilde{A}$ and the
variety $\widetilde{A}=A\cup \Sigma$ is smooth, compact and
connected. Indeed, let $\psi(t,p)=\{
z(t)=(z_{1}(t),\ldots,z_{5}(t)):t\in \mathbb{C},0< |t| <
\varepsilon \},$ be the orbit of the vector field (28) going
through the point $p\in A.$ Let $\Sigma _{p}\subset
\mathbb{P}^{15}$ be the surface element formed by the divisor
$\mathcal{D}$\ and the orbits going through $p,$ and set $\Sigma
\equiv \displaystyle{\cup_{p\in \mathcal{D}}\Sigma _{p}}.$\
Consider the Riemann surface $\mathcal{D}^{\prime
}=\mathcal{H}\cap \Sigma $\ where $\mathcal{H}\subset
\mathbb{P}^{15}$ is a hyperplane transversal to the direction of
the flow. If $\mathcal{D}^{\prime }$\ is smooth, then using the
implicit function theorem the surface $\Sigma $\ is smooth. But if
$\mathcal{D}^{\prime }$\ is singular at $0,$\ then $\Sigma $\
would be singular along the trajectory ($t-$axis) which go
immediately into the affine part A. Hence, A would be singular
which is a contradiction because A is the fibre of a morphism from
$\mathbb{C}^{5}$\ to $\mathbb{C}^{3}$\ and so smooth for almost
all the three constants of the motion $c_{k}.$ Next, let
$\overline{A}$ be the projective closure of A into
$\mathbb{P}^{5},$ let $Z=[Z_{0}:Z_{1}:\ldots:Z_{5}]\in
\mathbb{P}^{5}$ and let $I=\overline{A}\cap \{Z_0=0\}$ be the
locus at infinity. Consider the map $\overline{A}\subseteq
\mathbb{P}^{5} \rightarrow \mathbb{P}^{15},\text{ }Z\mapsto f(Z),$
where $f=(f_0,f_{1},...,f_{15}) \in \mathcal{L}(\mathcal{D})$ and
let $\widetilde{A}=f(\overline{A}).$ In a neighbourhood
$V(p)\subseteq \mathbb{P}^{15}$ of $p,$\ we have $\Sigma _{p}=
\widetilde{A}$ and $\Sigma _{p}\backslash \mathcal{D}\subseteq A.$
Otherwise there would exist an element of surface $\Sigma
_{p}^{\prime }\subseteq \widetilde{A}$ such that $\Sigma _{p}\cap
\Sigma _{p}^{\prime }=(t-axis)$, orbit $\psi
(t,p)=(t-axis)\backslash \ p\subseteq A,$ and hence A would be
singular along the $t-$axis which is impossible. Since the variety
$\overline{A}\cap \{Z_{0}\neq 0\}$ is irreducible and since the
generic hyperplane section $\mathcal{H}_{gen.}$\ of $\overline{A}$
is also irreducible, all hyperplane sections are connected and
hence I is also connected. Now, consider the graph $\Gamma
_{f}\subseteq \mathbb{P}^{5}\times \mathbb{P}^{15}$ of the map
$f,$ which is irreducible together with $\overline{A}.$ It follows
from the irreducibility of I that a generic hyperplane section
$\Gamma _{f}\cap \{\mathcal{H}_{gen.}\times \mathbb{P}^{15}\}$ is
irreducible, hence the special hyperplane section $\Gamma _{f}\cap
\{\{Z_{0}=0\}\times \mathbb{P}^{15}\}$ is connected and therefore
the projection map $proj_{\mathbb{P}^{15}}\{\Gamma _{f}\cap
\{\{Z_{0}=0\}\times \mathbb{P}^{15}\}\}=f(I)\equiv \mathcal{D},$
is connected. Hence, the variety $A\cup \Sigma =\widetilde{A}$ is
compact, connected and embeds smoothly into $\mathbb{P}^{15}$ via
$f.$ We wish to show that $\widetilde{A}$ is an abelian surface
equipped with two everywhere independent commuting vector fields.
For doing that, let $\phi^{\tau _{1}}$\ and $\phi^{\tau _{2}}$\ be
the flows corresponding to vector fields $X_{F_{1}}$ and
$X_{F_{2}}$. The latter are generated respectively by $F_1$ and
$F_2.$ For $p\in \mathcal{D}$ and for small $\varepsilon
> 0,$ $\phi^{\tau _{1}}(p), \forall \tau _{1},0<|\tau _{1}|<
\varepsilon ,$ is well defined and $\phi^{\tau _{1}}(p)\in
\widetilde{A}\backslash A.$ Then we may define $\phi^{\tau _{2}}$\
on $\widetilde{A}$ by $\phi^{\tau _{2}}(q)=\phi^{-\tau
_{1}}\phi^{\tau _{2}}\phi^{\tau _{1}}(q), q\in U(p)=\phi^{-\tau
_{1}}(U(\phi^{\tau _{1}}(p))),$ where $U(p)$ is a neighbourhood of
$p.$ By commutativity one can see that $\phi^{\tau _{2}}$ is
independent of $\tau _{1};$ $\phi^{-\tau _{1}-\varepsilon
_{1}}\phi^{\tau _{2}}\phi^{\tau _{1}+\varepsilon
_{1}}(q)=\phi^{-\tau _{1}}\phi^{-\varepsilon _{1}}\phi^{\tau
_{2}}\phi^{\tau _{1}}\phi^{\varepsilon _{1}}=\phi^{-\tau
_{1}}\phi^{\tau _{2}}\phi^{\tau _{1}}(q).$ We affirm that
$\phi^{\tau _{2}}(q)$ is holomorphic away from $\mathcal{D}.$ This
because $\phi^{\tau _{2}}\phi^{\tau _{1}}(q)$ is holomorphic away
from $\mathcal{D}$ and that $\phi^{\tau _{1}}$ is holomorphic in
$U(p)$ and maps bi-holomorphically $U(p)$ onto $U(\phi^{\tau
_{1}}(p)).$ Now, since the flows $\phi^{\tau _{1}}$ and
$\phi^{\tau _{2}}$ are holomorphic and independent on
$\mathcal{D},$ we can show along the same lines as in the
Arnold-Liouville theorem [15] that $\widetilde{A}$ is a complex
torus $\mathbb{C}^{2}/lattice$ and so in particular
$\widetilde{A}$ is a K\"{a}hler variety$.$ And that will done, by
considering the local diffeomorphism $\mathbb{C}^{2}\rightarrow
\widetilde{A}, (\tau _{1},\tau _{2})\mapsto \phi^{\tau
_{1}}\phi^{\tau _{2}}(p),$ for a fixed origin $p\in A.$ The
additive subgroup $\{(\tau _{1},\tau _{2})\in
\mathbb{C}^{2}:\phi^{\tau _{1}}\phi^{\tau _{2}}(p)=p\}$ is a
lattice of $\mathbb{C}^{2}$, hence
$\mathbb{C}^{2}/lattice\rightarrow \widetilde{A}$ is a
biholomorphic diffeomorphism and $\widetilde{A}$ is a K\"{a}hler
variety with K\"{a}hler metric given by $d\tau _{1}\otimes
d\overline{\tau }_{1}+d\tau _{2}\otimes d\overline{\tau }_{2}.$ As
mentioned in appendix A, a compact complex K\"{a}hler variety
having the required number as (its dimension) of independent
meromorphic functions is a projective variety. In fact, here we
have $\widetilde{A}\subseteq \mathbb{P}^{15}.$ Thus
$\widetilde{A}$ is both a projective variety and a complex torus
$\mathbb{C}^{2}/lattice$ and hence an abelian surface as a
consequence of Chow theorem. This completes the proof of the
theorem.

\begin{Rem}
a) Note that the reflection $\sigma$ on the affine variety $A$
amounts to the flip $\sigma :(z_1,z_2,z_3,z_4,z_5)\mapsto
(z_1,-z_2,z_3,-z_4,z_5),$ changing the direction of the commuting
vector fields. It can be extended to the (-Id)-involution about
the origin of $\mathbb{C}^2$ to the time flip $(t_1,t_2)\mapsto
(-t_1,-t_2)$ on $\widetilde{A}$, where $t_{1}$ and $t_{2}$ are the
time coordinates of each of the flows $X_{{F}_{1}}$ and
$X_{{F}_{2}}.$ The involution $\sigma $ acts on the parameters of
the Laurent solution (30) as follows $\sigma
:(t,\alpha,\beta,\gamma,\theta)\longmapsto
(-t,-\alpha,-\beta,-\gamma,\theta),$ interchanges the Riemann
surfaces $\mathcal{C}_{\varepsilon}$ and the linear space
$\mathcal{L}$ can be split into a direct sum of even and odd
functions. Geometrically, this involution interchanges
$\mathcal{C}_{1}$ and $\mathcal{C}_{-1},$ i.e.,
$\mathcal{C}_{-1}=\sigma \mathcal{C}_{1}.$\\
b) Consider on $\widetilde{A}$ the holomorphic 1-forms $dt_1$ and
$dt_2$ defined by $dt_i(X_{F_j})=\delta_{ij},$ where $X_{F_1}$ and
$X_{F_2}$ are the vector fields generated respectively by $F_1$
and $F_2.$ Taking the differentials of $\zeta=1/z_1$ and
$\xi=z_1/z_2$ viewed as functions of $t_1$ and $t_2,$ using the
vector fields and the Laurent series (31) and solving linearly for
$dt_1$ and $dt_2,$ we obtain the holomorphic differentials
\begin{eqnarray}
\omega_1&=&dt_{1}|_{\mathcal{C}_{\varepsilon}}=\frac{1}{\triangle}(\frac{\partial
\xi}{\partial t_2}d\zeta-\frac{\partial \zeta}{\partial
t_2}d\xi)|_{\mathcal{C}_{\varepsilon}} =
\frac{8}{\alpha \left(-4\beta +\alpha ^{3}\right) }d\alpha ,\nonumber\\
\omega_2&=&dt_{2}|_{\mathcal{C}_{\varepsilon}}=\frac{1}{\triangle}(\frac{-\partial
\xi}{\partial t_1}d\zeta-\frac{\partial \zeta}{\partial
t_1}d\xi)|_{\mathcal{C}_{\varepsilon}}=\frac{2}{\left(-4\beta
+\alpha ^{3}\right) ^{2}}d\alpha ,\nonumber
\end{eqnarray}
with $\Delta \equiv {\frac{\partial \zeta}{\partial
t_1}\frac{\partial \xi}{\partial t_2}-\frac{\partial
\zeta}{\partial t_2}\frac{\partial \xi}{\partial t_1}}.$ The
zeroes of $\omega_2$ provide the points of tangency of the vector
field $X_{F_{1}}$ to $\mathcal{C}_\varepsilon.$ We have
$\frac{\omega_1}{\omega_2}=\frac{4}{\alpha }\left( -4\beta +\alpha
^{3}\right),$ and $X_{F_{1}}$ is tangent to
$\mathcal{H}_\varepsilon$ at the point covering $\alpha =\infty.$
\end{Rem}

\subsection{The Hénon-Heiles system}
The Hénon-Heiles system
\begin{equation}\label{eqn:euler}
\dot q_{1}=\frac{\partial H}{\partial p_{1}},\quad
q_{2}=\frac{\partial H}{\partial p_{2}},\quad \dot
p_{1}=-\frac{\partial H}{\partial q_{1}},\quad \dot
p_{2}=-\frac{\partial H}{\partial q_{2}},
\end{equation}
with $$H\equiv H_{1}=\frac{1}{2}\left(
p_{1}^{2}+p_{2}^{2}+aq_{1}^{2}+bq_{2}^{2}\right)
+q_{1}^{2}q_{2}+6q_{2}^{3},$$ has another constant of motion
$$H_{2}=q_{1}^{4}+4q_{1}^{2}q_{2}^{2}-4p_{1}\left( p_{1}q_{2}-p_{2}q_{1}\right)
+4aq_{1}^{2}q_{2}+\left( 4a-b\right) \left(
p_{1}^{2}+aq_{1}^{2}\right),$$ where $a,$ $b,$ are constant
parameters and $q_{1},q_{2},p_{1},p_{2}$ are canonical coordinates
and momenta, respectively. First studied as a mathematical model
to describe the chaotic motion of a test star in an axisymmetric
galactic mean gravitational field this system is widely explored
in other branches of physics. It well-known from applications in
stellar dynamics, statistical mechanics and quantum mechanics. It
provides a model for the oscillations of atoms in a three-atomic
molecule. The system $\left(34\right) $ possesses Laurent series
solutions depending on $3$ free parameters $\alpha ,\beta ,\gamma
,$ namely
\begin{eqnarray}
q_{1}&=&\frac{\alpha }{t}+\left( \frac{\alpha
^{3}}{12}+\allowbreak \frac{\alpha A}{2}-\frac{\alpha
B}{12}\right) t+\beta t^{2}+q_{1}^{\left( 4\right)
}t^{3}+q_{1}^{\left( 5\right)
}t^{4}+q_{1}^{\left( 6\right) }t^{5}+\cdots,\nonumber\\
q_{2}&=&-\frac{1}{t^{2}}+\frac{\alpha
^{2}}{12}-\frac{B}{12}+\left( \frac{\alpha ^{4}}{48}+\frac{\alpha
^{2}A}{10}-\frac{\alpha ^{2}B}{60}-\frac{B^{2}}{240}\right)
t^{2}+\frac{\alpha \beta }{3}t^{3}+\gamma t^{4}+\cdots,\nonumber
\end{eqnarray}
where $p_{1}=\overset{.}{q}_{1},\text{ }p_{2}=\overset{.}{q}_{2}$
and
\begin{eqnarray}
q_{1}^{\left( 4\right) }&=&\frac{\alpha AB}{24}-\frac{\alpha
^{5}}{72}+\frac{11\alpha ^{3}B}{720}-\frac{11\alpha
^{3}A}{120}-\frac{\alpha B^{2}}{720}-\frac{\alpha
A^{2}}{8},\nonumber\\
q_{1}^{\left( 5\right) }&=&-\frac{\beta \alpha
^{2}}{12}+\frac{\beta B}{60}-\frac{A\beta }{10},\nonumber\\
q_{1}^{\left( 6\right) }&=&-\frac{\alpha \gamma }{9}-\frac{\alpha
^{7}}{15552}-\frac{\alpha ^{5}A}{2160}+\frac{\alpha
^{5}B}{12960}+\frac{\alpha ^{3}B^{2}}{25920}+\frac{\alpha
^{3}A^{2}}{1440}-\frac{\alpha ^{3}AB}{4320}+\frac{\alpha
AB^{2}}{1440}\nonumber\\
&&-\frac{\alpha B^{3}}{19440}-\frac{\alpha
A^{2}B}{288}+\frac{\alpha A^{3}}{144}.\nonumber
\end{eqnarray}
Let $\mathcal{D}$ be the pole solutions restricted to the surface
$$M_{c}=\overset{2}{\underset{i=1}{\bigcap }}\left\{ x\equiv
(q_1,q_2,p_1,p_2)\in \mathbb{C}^{4}, H_{i}\left( x\right)
=c_{i}\right\},$$ to be precise $\mathcal{D}$ is the closure of
the continuous components of the set of Laurent series solutions
$x\left( t\right)$ such that $H_{i}\left( x\left( t\right) \right)
=c_{i},\text{ }1\leq i\leq 2$, i.e., $
\mathcal{D}=t^{0}-\text{coefficient of }M_c. $ Thus we find an
algebraic curve defined by
\begin{equation}\label{eqn:euler}
\mathcal{D} : \beta ^{2}=P_{8}(\alpha),
\end{equation}
where $$
P_{8}\left( \alpha \right) =-\frac{7}{15552}\alpha ^{8}-\frac{1}{432}%
\left( 5A-\frac{13}{18}B\right) \alpha ^{6}
-\frac{1}{36}\left( \frac{671}{15120}B^{2}+\frac{17}{7}A^{2}-\frac{943}{%
1260}BA\right) \alpha ^{4} $$
$$-\frac{1}{36}\left( 4A^{3}-\frac{1}{2520}B^{3}-\frac{13}{6}A^{2}B+\frac{2}{%
9}AB^{2}-\frac{10}{7}c_{1}\right) \alpha ^{2}+\frac{1}{36}c_{2}.
$$ The curve $\mathcal{D}$ determined by an eight-order equation
is smooth, hyperelliptic and its genus is $3$. Moreover, the map
\begin{equation}\label{eqn:euler}
\sigma :\mathcal{D}\longrightarrow \mathcal{D},\text{ }(\beta
,\alpha ) \longmapsto (\beta ,-\alpha ),
\end{equation}
is an involution on $\mathcal{D}$ and the quotient
$\mathcal{E}=\mathcal{D}/\sigma $ is an elliptic curve defined by
\begin{equation}\label{eqn:euler}
\mathcal{E}: \beta ^{2}=P_{4}(\zeta),
\end{equation}
where $P_{4}\left( \zeta \right) $ is the degree $4$ polynomial in
$\zeta =\alpha ^{2}$ obtained from $\left(35\right) .$ The
hyperelliptic curve $\mathcal{D}$ is thus a $2$-sheeted ramified
covering of the elliptic curve $\mathcal{E}\left(37\right) ,$
\begin{equation}\label{eqn:euler}
\rho: \mathcal{D}\longrightarrow \mathcal{E},\text{ }(\beta
,\alpha )\longmapsto (\beta ,\zeta ),
\end{equation}
ramified at the four points covering $\zeta =0$ and $\infty .$ The
affine surface $M_c$ completes into an abelian surface
$\widetilde{M}_c$, by adjoining the divisor $\mathcal{D}$. The
latter defines on $\widetilde{M}_c$ a polarization $(1,2).$ The
divisor $2\mathcal{D}$ is very ample and the functions $1,\text{
}y_{1},\text{ }y_{1}^{2},\text{ }y_{2},\text{ }x_{1}, \text{
}x_{1}^{2}+y_{1}^{2}y_{2},\text{ }x_{2}y_{1}-2x_{1}y_{2},\text{ }
x_{1}x_{2}+2Ay_{1}y_{2}+2y_{1}y_{2}^{2},$ embed $\widetilde{M}_c$
smoothly into $\mathbb{CP}^7$ with polarization $(2,4).$ Then the
system (34) is algebraically completely integrable and the
corresponding flow evolues on an abelian surface $\widetilde{M}_c
= \mathbb{C}^2 /\mbox{lattice},$ where the lattice is generated by
the period matrix $\left(
\begin{array}{llll}
2 & 0 & a & c \\
0 & 4 & c & b
\end{array}
\right)$, $\text{ Im}\left(
\begin{array}{ll}
a & c \\
c & b
\end{array}
\right) >0.$

\begin{Theo}
The abelian surface $\widetilde{M}_{c}$ which completes the affine
surface $M_{c}$ is the dual Prym variety $Prym^{*}\left(
\mathcal{D}/\mathcal{E}\right) $ of the genus $3$ hyperelliptic
curve $\mathcal{D}$ (35) for the involution $\sigma$ interchanging
the sheets of the double covering $\rho $ (38) and the problem
linearizes on this variety.
\end{Theo}
\emph{Proof}. Let $\left(
a_{1},a_{2},a_{3},b_{1},b_{2},b_{3}\right) $ be a canonical
homology basis of $\mathcal{D}$ such that $\sigma \left(
a_{1}\right) =a_{3},$ $\sigma \left( b_{1}\right) =b_{3},$ $\sigma
\left( a_{2}\right) =-a_{2},$ $\sigma \left( b_{2}\right)
=-b_{2},$ for the involution $\sigma $ (36). As a basis of
holomorphic differentials $\omega _{0},\omega _{1},\omega _{2}$ on
the curve $\mathcal{D}$ (35) we take the differentials $\omega
_{1}=\frac{\alpha ^{2}d\alpha }{\beta },\omega _{2}= \frac{d\alpha
}{\beta }, \omega _{3}=\frac{\alpha d\alpha }{\beta},$ and
obviously $\sigma ^{*}(\omega _{1})=-\omega _{1},\sigma
^{*}(\omega _{2})=-\omega _{2},.\sigma ^{*}(\omega _{3})=\omega
_{3}.$ Recall that the Prym variety $Prym\left(
\mathcal{D}/\mathcal{E}\right) $\ \ is a subabelian variety of the
Jacobi variety
$Jac(\mathcal{D})=Pic^{0}(\mathcal{D})=H^{1}(\mathcal{O}_{\mathcal{D}})\
/\ H^{1}(\mathcal{D},\mathbb{Z})$ constructed from the double
cover $\rho $ $:$\ the involution $\sigma $\ \ on $\mathcal{D}$\
interchanging sheets, extends by linearity to a map $\sigma
:Jac(\mathcal{D})\rightarrow Jac(\mathcal{D})$ and up to some
points of order two, $Jac(\mathcal{D})$ splits into an even part
and an odd part : the even part is an elliptic curve (the quotient
of $\mathcal{D}$ by $\sigma $, i.e., $\mathcal{E}$ (18)) and the
odd part is a $2-$dimensional abelian surface $Prym\left(
\mathcal{D}/\mathcal{E}\right).$ We consider the period matrix
$\Omega$ of $Jac(D)$
$$\Omega=\left(
\begin{array}{cccccc}
\int_{a_{1}}\omega _{1}&\int_{a_{2}}\omega _{1}&\int_{a_{3}}\omega
_{1}&\int_{b_{1}}\omega _{1}&\int_{b_{2}}\omega _{1}&\int_{b_{3}}\omega _{1}\\
\int_{a_{1}}\omega _{2}&\int_{a_{2}}\omega _{2}&\int_{a_{3}}\omega
_{2}&\int_{b_{1}}\omega _{2}&\int_{b_{2}}\omega _{2}&\int_{b_{3}}\omega _{2}\\
\int_{a_{1}}\omega _{3}&\int_{a_{2}}\omega _{3}&\int_{a_{3}}\omega
_{3}&\int_{b_{1}}\omega _{3}&\int_{b_{2}}\omega
_{3}&\int_{b_{3}}\omega _{3}
\end{array}
\right).
$$
Then,
$$\Omega=\left(
\begin{array}{cccccc}
\int_{a_{1}}\omega _{1}&\int_{a_{2}}\omega
_{1}&-\int_{a_{1}}\omega _{1}
&\int_{b_{1}}\omega _{1}&\int_{b_{2}}\omega _{1}&-\int_{b_{1}}\omega _{1}\\
\int_{a_{1}}\omega _{2}&\int_{a_{2}}\omega
_{2}&-\int_{a_{1}}\omega _{2}
&\int_{b_{1}}\omega _{2}&\int_{b_{2}}\omega _{2}&-\int_{b_{1}}\omega _{2}\\
\int_{a_{1}}\omega _{3}&0&\int_{a_{1}}\omega
_{3}&\int_{b_{1}}\omega _{3}&0&\int_{b_{1}}\omega _{3}
\end{array}
\right),
$$
and therefore the period matrices of $Jac(\mathcal{E})$(i.e.,
$\mathcal{E}$), $Prym(\mathcal{D}/\mathcal{E})$ and
$Prym^*(\mathcal{D}/\mathcal{E})$ are respectively
$\Delta=(\int_{a_{1}}\omega _{3}\quad \int_{b_{1}}\omega _{3}),$
$$
\Gamma=\left(
\begin{array}{cccc}
2\int_{a_{1}}\omega _{1}&\int_{a_{2}}\omega
_{1}&2\int_{b_{1}}\omega _{1}
&\int_{b_{2}}\omega _{1}\\
2\int_{a_{1}}\omega _{2}&\int_{a_{2}}\omega
_{2}&2\int_{b_{1}}\omega _{2} &\int_{b_{2}}\omega _{2}
\end{array}
\right),
$$ and
$$
\Gamma^*=\left(
\begin{array}{cccc}
\int_{a_{1}}\omega _{1}&\int_{a_{2}}\omega _{1}&\int_{b_{1}}\omega
_{1}
&\int_{b_{2}}\omega _{1}\\
\int_{a_{1}}\omega _{2}&\int_{a_{2}}\omega _{2}&\int_{b_{1}}\omega
_{2} &\int_{b_{2}}\omega _{2}
\end{array}
\right).
$$
Let $L_{\Omega
}=\{\sum_{i=1}^{3}m_{i}\int_{a_{i}}\left(\begin{array}{c}
\omega_1\\\omega_2 \\\omega_3
\end{array}\right )
+n_{i}\int_{b_{i}}\left(\begin{array}{c} \omega_1\\\omega_2
\\\omega_3
\end{array}\right) :m_{i},n_{i}\in \mathbb{Z}\} ,$
be the period lattice associated to $\Omega.$ Let us denote also
by $L_{\Delta},$ the period lattice associated $\Delta.$ We have
the following diagram
$$
\begin{array}{ccccccccc}
&&&&0&&&\\
&&&&\downarrow&&&&\\
&&&&\mathcal{E}&&\mathcal{D}&&\\
&&&&\quad \downarrow \varphi^*&\swarrow&\quad \downarrow \varphi&&\\
0&\longrightarrow&\ker
N_\varphi&\longrightarrow&Prym(\mathcal{D}/\mathcal{E} \oplus
\mathcal{E}=Jac(\mathcal{D})&\overset{N_{\varphi
}}{\longrightarrow }&
\mathcal{E}&\longrightarrow&0\\
&&&\searrow \tau&\downarrow&&&&\\
&&&&\widetilde{M}_c=M_c\cup 2\mathcal{D}\simeq \mathbb{C}^2/\mbox{lattice}&&&&\\
&&&&\downarrow&&&&\\
&&&&0&&&&
\end{array}
$$
The polarization map $\tau
:Prym(\mathcal{D}/\mathcal{E})\longrightarrow
\widetilde{M}_{c}=Prym^{*}(\mathcal{D}/\mathcal{E}),$ has kernel
$(\varphi ^{*}\mathcal{E})\simeq \mathbb{Z}_{2}\times
\mathbb{Z}_{2}$ and the induced polarization on
$Prym(\mathcal{D}/\mathcal{E})$ is of type (1,2). Let
$\widetilde{M}_{c}\rightarrow
\mathbb{C}^{2}/L_{\Lambda}:p\curvearrowright
\int_{p_{0}}^{p}\binom{dt_{1}}{dt_{2}},$ be the uniformizing map
where $dt_1, dt_2$ are two differentials on $\widetilde{M}_c$
corresponding to the flows generated respectively by $H_1, H_2$
such that : $dt_1|_\mathcal{D}=\omega_1$ and
$dt_2|_\mathcal{D}=\omega_2$,
$$L_{\Lambda
}=\{\sum_{k=1}^{4}n_{k}\left(\begin{array}{c}
\int_{\nu_{k}}dt_1\\\int_{\nu_{k}}dt_2
\end{array}\right ):n_{k}\in \mathbb{Z}\} ,$$
is the lattice associated to the period matrix
$$
\Lambda=\left(
\begin{array}{cccc}
\int_{\nu_{1}}dt _{1}&\int_{\nu_{2}}dt _{1}&\int_{\nu_{4}}dt _{1}
&\int_{\nu_{4}}dt _{1}\\
\int_{\nu_{1}}dt _{2}&\int_{\nu_{2}}dt _{2}&\int_{\nu_{3}}dt _{2}
&\int_{\nu_{4}}dt _{2}
\end{array}
\right),
$$
and $(\nu_{1},\nu_{2},\nu_{3},\nu_{4})$ is a basis of
$H_{1}(\widetilde{M}_{c},\mathbb{Z})$. By the Lefschetz theorem on
hyperplane section [9], the map
$H_{1}(\mathcal{D},\mathbb{Z})\longrightarrow
H_{1}(\widetilde{M}_{c},\mathbb{Z})$\ induced by the inclusion
$\mathcal{D}\hookrightarrow $ $\widetilde{M}_{c}$\ is surjective
and consequently we can find $4$\ cycles $\nu _{1},\nu _{2},\nu
_{3},\nu _{4}$\ on the curve $\mathcal{D}$ such that
$$
\Lambda=\left(
\begin{array}{cccc}
\int_{\nu_{1}}\omega _{1}&\int_{\nu_{2}}\omega
_{1}&\int_{\nu_{4}}\omega _{1}
&\int_{\nu_{4}}\omega _{1}\\
\int_{\nu_{1}}\omega _{2}&\int_{\nu_{2}}\omega
_{2}&\int_{\nu_{3}}\omega _{2} &\int_{\nu_{4}}\omega _{2}
\end{array}
\right),
$$
and $L_{\Lambda }=\{\sum_{k=1}^{4}n_{k}\left(\begin{array}{c}
\int_{\nu_{k}}\omega_1\\\int_{\nu_{k}}\omega_2
\end{array}\right ):n_{k}\in \mathbb{Z}\}.$
The  cycles $\nu _{1},\nu _{2},\nu _{3},\nu _{4}$ in $\mathcal{D}$
which we look for are $a_{1},b_{1},a_{2},b_{2}$ and they generate
$H_{1}(\widetilde{M}_{c},\mathbb{Z})$ such that
$$
\Lambda=\left(
\begin{array}{cccc}
\int_{a_{1}}\omega _{1}&\int_{b_{1}}\omega _{1}&\int_{a_{2}}\omega
_{1}&\int_{b_{2}}\omega _{1}\\
\int_{a_{1}}\omega _{2}&\int_{b_{1}}\omega _{2}&\int_{a_{2}}\omega
_{2} &\int_{b_{2}}\omega _{2}
\end{array}
\right),
$$
is a Riemann matrix. We show that $\Lambda=\Gamma^*$ ,i.e., the
period matrix of $Prym^*(\mathcal{D}/\mathcal{E})$ dual of
$Prym(\mathcal{D}/\mathcal{E})$. Consequently $\widetilde{M}_{c}$
and $Prym^*(\mathcal{D}/\mathcal{E})$ are two abelian varieties
analytically isomorphic to the same complex torus
$\mathbb{C}^{2}/L_{\Lambda}.$ By Chow's theorem [9],
$\widetilde{\mathcal{A}}_{c}$\ and
$Prym^*(\mathcal{D}/\mathcal{E})$ are then algebraically
isomorphic.

\subsection{The Kowalewski rigid body motion}
The motion for the Kowalewski's top is governed by the equations
\begin{equation}\label{euler:eqn}
\overset{.}{m}=m\wedge \lambda m+\gamma \wedge l,\qquad
\overset{.}{\gamma }=\gamma \wedge \lambda m,
\end{equation}
where $m,\gamma $ and $l$ denote respectively the angular
momentum, the directional cosine of the $z$-axis (fixed in space),
the center of gravity which after some rescaling and normalization
may be taken as $l=\left( 1,0,0\right) $ and $\lambda m=\left(
m_{1}/2,m_{2}/2,m_{3}/2\right) .$ The system (39) can be written
\begin{eqnarray}
\overset{.}{m}_{1}&=&m_{2}\text{ }m_{3},\qquad \qquad \qquad
\overset{.}{\gamma
}_{1}=2\text{ }m_{3}\gamma _{2}-m_{2}\gamma _{3},\nonumber\\
\overset{.}{m}_{2}&=&-\text{ }m_{1}\text{ }m_{3}+2\gamma
_{3},\qquad \overset{.}{\gamma }_{2}=\text{ }m_{1}\gamma
_{3}-2m_{3}\gamma _{1},\\
\overset{.}{m}_{3}&=&-2\gamma_{2}\text{,}\qquad \qquad \qquad
\quad \overset{.}{\gamma }_{3}=m_{2}\gamma _{1}-m_{1}\gamma
_{2},\nonumber
\end{eqnarray}
with constants of motion
\begin{eqnarray}
H_{1}&=&\frac{1}{2}\left( m_{1}^{2}+m_{2}^{2}\right) +m_{3}^{2}+2\gamma _{1}=c_{1},\nonumber\\
H_{2}&=&m_{1}\gamma _{1}+m_{2}\gamma _{2}+m_{3}\gamma _{3}=c_{2},\\
H_{3}&=&\gamma _{1}^{2}+\gamma _{2}^{2}+\gamma _{3}^{2}=c_{3}=1,\nonumber\\
H_{4}&=&\left( \left( \frac{m_{1}+im_{2}}{2}\right) ^{2}-\left(
\gamma _{1}+i\gamma _{2}\right) \right) \left( \left(
\frac{m_{1}-im_{2}}{2}\right) ^{2}-\left( \gamma _{1}-i\gamma
_{2}\right) \right) =c_{4}.\nonumber
\end{eqnarray}
The system (40) admits two distinct families of Laurent series
solutions :
$$
m_{1}\left( t\right)=\left\{\begin{array}{rl}
&\frac{\alpha_{1}}{t}+i\left( \alpha _{1}^{2}-2\right) \alpha
_{2}+ \circ\left( t\right),\\
&\frac{\alpha _{1}}{t}-i\left( \alpha _{1}^{2}-2\right) \alpha
_{2}+ \circ\left( t\right),
\end{array}\right.
\quad \gamma _{1}\left( t\right)= \left\{\begin{array}{rl}
&\frac{1}{2t^{2}}+\circ\left( t\right) ,\\
&\frac{1}{2t^{2}}+\circ\left( t\right) ,\end{array}\right.$$

$$ m_{2}\left( t\right)= \left\{\begin{array}{rl} &\frac{i\alpha
_{1}}{t}-\alpha _{1}^{2}\alpha _{2}+\circ\left( t\right) ,\\
&\frac{-i\alpha _{1}}{t}-\alpha _{1}^{2}\alpha _{2} +\circ\left(
t\right) ,\end{array}\right. \qquad \gamma _{2}\left( t\right)=
\left\{\begin{array}{rl}
&\frac{i}{2t^{2}}+\circ\left( t\right) ,\\
&\frac{-i}{2t^{2}}+\circ\left( t\right) ,\end{array}\right.$$

$$m_{3}\left( t\right)= \left\{\begin{array}{rl}
&\frac{i}{t}+\alpha _{1}\alpha _{2}
+\circ\left( t\right) ,\\
&\frac{-i}{t}+\alpha _{1}\alpha _{2}+\circ\left( t\right) ,
\end{array}\right.\qquad
\gamma _{3}\left( t\right) = \left\{\begin{array}{rl}
&\frac{\alpha _{2}}{t}+\circ\left( t\right) ,\\
&\frac{\alpha _{2}}{t}+\circ\left( t\right) ,\end{array}\right.$$
which depend on $5$ free parameters $\alpha _{1},...,$ $\alpha
_{5.}$ By substituting these series in the constants of the motion
$H_{i}$ (41), one eliminates three parameters linearly, leading to
algebraic relation between the two remaining parameters, which is
nothing but the equation of the divisor $\mathcal{D}$ along which
the $m_{i},\gamma _{i}$ blow up. Since the system (40) admits two
families of Laurent solutions, then $\mathcal{D}$ is a set of two
isomorphic curves of genus $3,$
$\mathcal{D}=\mathcal{D}_{1}+\mathcal{D}_{-1}:$
\begin{equation}\label{eqn:euler}
\mathcal{D}_{\varepsilon }:\text{ }P\left( \alpha _{1},\alpha
_{2}\right) =\left( \alpha _{1}^{2}-1\right) \left( \left( \alpha
_{1}^{2}-1\right) \alpha _{2}^{2}-P\left( \alpha _{2}\right)
\right)+c_{4}=0,
\end{equation}
where $P\left( \alpha _{2}\right) =c_{1}\alpha
_{2}^{2}-2\varepsilon c_{2}\alpha _{2}-1$ and $\varepsilon =\pm
1.$ Each of the curve $\mathcal{D}_{\varepsilon }$ is a $2-1$
ramified cover $\left( \alpha _{1},\alpha _{2},\beta \right) $ of
elliptic curves $\mathcal{D}_{\varepsilon }^{0}:$
\begin{equation}\label{eqn:euler}
\mathcal{D}_{\varepsilon }^{0}:\beta ^{2}=P^{2}\left( \alpha
_{2}\right) -4c_{4}\alpha _{2}^{4},
\end{equation}
ramified at the $4$ points $\alpha _{1}=0$ covering the $4$ roots
of $P\left( \alpha _{2}\right) =0.$ It was shown $\left[12\right]
$ that each divisor $\mathcal{D}_{\varepsilon }$ is ample and
defines a polarization $\left( 1,2\right) ,$ whereas the divisor
$\mathcal{D},$ of geometric genus $9,$ is very ample and defines a
polarization $\left( 2,4\right)$. The affine surface
$M_{c}=\bigcap_{i=1}^{4}\left\{ H_{i}=c_{i}\right\} \subset
\mathbb{C}^{6},$ defined by putting the four invariants (41) of
the Kowalewski flow (40) equal to generic constants, is the affine
part of an abelian surface $\widetilde{M_{c}}$ with
\begin{eqnarray}
\widetilde{M_{c}}\text{ }\backslash \text{
}M_{c}=\mathcal{D}&=&\text{one
genus 9 curve consisting of two genus 3 }\nonumber\\
&&\text{curves }\mathcal{D}_{\varepsilon }\text{ }(42) \text{
intersecting in 4 points. Each }\nonumber\\
&&\mathcal{D}_{\varepsilon }\text{ is a double cover of an
elliptic curve
}\mathcal{D}_{\varepsilon }^{0}\text{ }(43) \nonumber\\
&&\text{ramified at 4 points.}\nonumber
\end{eqnarray}
Moreover, the Hamiltonian flows generated by the vector fields
$X_{H_{1}}$ and $X_{H_{4}}$ are straight lines on
$\widetilde{M_{c}}.$ The $8$ functions $1$ $,$ $f_{1}=m_{1}$ $,$
$f_{2}=m_{2}$ $,$ $f_{3}=m_{3}$ $,$ $f_{4}=\gamma _{3}$ $,$
$f_{5}=f_{1}^{2}+f_{2}^{2}$ $,$  $f_{6}=4f_{1}f_{4}-f_{3}f_{5}$
$,$ $f_{7}=\left( f_{2}\gamma _{1}-f_{1}\gamma _{2}\right)
f_{3}+2f_{4}\gamma _{2},$ form a basis of the vector space of
meromorphic functions on $\widetilde{M_{c}}$ with at worst a
simple pole along $\mathcal{D}$ Moreover, the map
$$\widetilde{M_{c}}\simeq \mathbb{C}^{2}/Lattice\rightarrow
\mathbb{CP}^{7}\text{ },\text{ }\left( t_{1},t_{2}\right) \mapsto
\left[ \left( 1,f_{1}\left( t_{1},t_{2}\right) ,...,f_{7}\left(
t_{1},t_{2}\right) \right) \right] ,$$ is an embedding of
$\widetilde{M_{c}}$ into $\mathbb{CP}^{7}.$ Following the method
(theorem 13), we obtain the following theorem :
\begin{Theo}
The tori $\widetilde{M_{c}}$ can be identified as
$\widetilde{M_{c}}=Prym^*(\mathcal{D}_\varepsilon/\mathcal{D}_\varepsilon
^0$, i.e., dual of $Prym(\mathcal{D}_{\varepsilon
}/\mathcal{D}_{\varepsilon }^{0})$ and the problem linearizes on
this Prym variety.
\end{Theo}

\subsection{Kirchhoff's equations of motion of a solid in an ideal fluid}
The Kirchhoff's equations of motion of a solid in an ideal fluid
have the form
\begin{eqnarray}
\dot p_{1}&=&p_{2}\frac{\partial H}{\partial l_{3}}-p_{3}
\frac{\partial H}{\partial l_{2}}\text{ ,\qquad }\dot l_{1}=
p_{2}\frac{\partial H}{\partial p_{3}}-p_{3}\frac{\partial
H}{\partial p_{2}}+l_{2}\frac{\partial H}{\partial
l_{3}}-l_{3}\frac{\partial H}{\partial l_{2}},\nonumber\\
\dot p_{2}&=&p_{3}\frac{\partial H}{\partial l_{1}}-p_{1}
\frac{\partial H}{\partial l_{3}}\text{ ,\qquad }\dot l_{2}=
p_{3}\frac{\partial H}{\partial p_{1}}-p_{1}\frac{\partial
H}{\partial p_{3}}+ l_{3}\frac{\partial H}{\partial
l_{1}}-l_{1}\frac{\partial H}{\partial l_{3}},\\
\dot p_{3}&=&p_{1}\frac{\partial H}{\partial l_{2}}-p_{2}
\frac{\partial H}{\partial l_{1}}\text{ ,\qquad }\dot l_{3}=
p_{1}\frac{\partial H}{\partial p_{2}}-p_{2}\frac{\partial
H}{\partial p_{1}}+ l_{1}\frac{\partial H}{\partial
l_{2}}-l_{2}\frac{\partial H}{\partial l_{1}},\nonumber
\end{eqnarray}
where $(p_{1},p_{2},p_{3})$\ is the velocity of a point fixed
relatively to the solid, $(l_{1},l_{2},l_{3})$\ the angular
velocity of the body expressed with regard to a frame of reference
also fixed relatively to the solid and $H$\ is the hamiltonian.
These equations \ can be regarded as the equations of the
geodesics of the right-invariant metric on the group $E\left(
3\right) =SO\left( 3\right) \times \mathbb{R}^{3}$\ of motions of
3-dimensional euclidean space $\mathbb{R}^{3},$\ generated by
rotations and translations. Hence the motion has the trivial
coadjoint orbit invariants $\langle p,p\rangle$ and $\langle
p,l\rangle.$ As it turns out, this is a special case of a more
general system of equations written as
$$\dot x=x\wedge \frac{\partial H}{\partial x}
+y\wedge \frac{\partial H}{\partial y},\quad \dot y=y\wedge
\frac{\partial H}{\partial x} +x\wedge \frac{\partial H}{\partial
y},$$ where $x=\left( x_{1},x_{2},x_{3}\right) \in \mathbb{R}^{3}$
et $y=\left( y_{1},y_{2},y_{3}\right) \in \mathbb{R}^{3}.$ The
first set can be obtained from the second by putting
$(x,y)=(l,p/\varepsilon)$ and letting $\varepsilon\rightarrow 0.$
The latter set of equations is the geodesic flow on $SO(4)$ for a
left invariant metric defined by the quadratic form $H.$ In
Clebsch's case, equations (44) have the four invariants :
\begin{eqnarray}
H_{1}&=&H=\frac{1}{2}\left(
a_{1}p_{1}^{2}+a_{2}p_{2}^{2}+a_{3}p_{3}^{2}+b_{1}l_{1}^{2}+
b_{2}l_{2}^{2}+b_{3}l_{3}^{2}\right),\nonumber\\
H_{2}&=&p_{1}^{2}+p_{2}^{2}+p_{3}^{2},\nonumber\\
H_{3}&=&p_{1}l_{1}+p_{2}l_{2}+p_{3}l_{3},\nonumber\\
H_{4}&=&\frac{1}{2}\left(
b_{1}p_{1}^{2}+b_{2}p_{2}^{2}+b_{3}p_{3}^{2}+\varrho \left(
l_{1}^{2}+l_{2}^{2}+l_{3}^{2}\right) \right) ,\nonumber
\end{eqnarray}
with $\frac{a_{2}-a_{3}}{b_{1}}+\frac{a_{3}-a_{1}}{b_{2}}+
\frac{a_{1}-a_{2}}{b_{3}}=0,$ and the constant $\varrho $\
satisfies the conditions $\varrho =\frac{b_{1}\left(
b_{2}-b_{3}\right) }{a_{2}-a_{3}}=\frac{b_{2}\left(
b_{3}-b_{1}\right) }{a_{3}-a_{1}}=\frac{b_{3}\left(
b_{1}-b_{2}\right) }{a_{1}-a_{2}}.$ The system (44) can be written
in the form (11) with $m=6;$\ to be precise
\begin{equation}\label{eqn:euler}
\dot{x}=f\left( x\right) \equiv J\frac{\partial H}{\partial
x}\text{ },\text{
}x=(p_{1},p_{2},p_{3},l_{1},l_{2},l_{3})^{\intercal },
\end{equation}
where $$J=\left(\begin{array}{cc}
O&P\\
P&L
\end{array}\right),
P=\left(\begin{array}{ccc}
0&-p_{3}&p_{2}\\
p_{3}&0&-p_{1}\\
-p_{2}&p_{1}&0
\end{array}\right),
L=\left(\begin{array}{ccc}
0&-l_{3}&l_{2}\\
l_{3}&0&-l_{1}\\
-l_{2}&p_{1}&0
\end{array}\right).
$$
Consider points at infinity which are limit points of trajectories
of the flow. In fact, there is a Laurent decomposition of such
asymptotic solutions,
\begin{equation}\label{eqn:euler}
x\left( t\right) =t^{-1}\left( x^{\left( 0\right) }+ x^{\left(
1\right) }t+x^{\left( 2\right) }t^{2}+...\right),
\end{equation}
which depend on $\dim (phase$\ $space)-1=5$\ free parameters$.$\
Putting (46) into (45), solving inductively for the $x^{\left(
k\right)},$ one finds at the $0^{th}$\ step a non-linear equation,
$x^{\left( 0\right) }+f(x^{\left( 0\right) })=0,$ and at the
$k^{th}$\ step, a linear system of equations,
$$
(L-kI) x^{(k)}= \left\{
\begin{array}{rl}
0 & \mbox{for}\quad k=1\\
\text{quadratic polynomial in }x^{\left( 1\right) },...,x^{\left(
k\right) }& \mbox{for}\quad k \geq 1,
\end{array}
\right.
$$
where $L$\ denotes the jacobian map of the non-linear equation
above.\ One parameter appear at the $0^{th}$\ step, i.e., in the
resolution of the non-linear equation\ and the $4$\ remaining ones
at the $k^{th}$\ step, $k=1,...,4.$ Taking into account only
solutions trajectories lying on the invariant surface
$M_{c}=\overset{4}{\underset{i=1}{\bigcap }}\left\{ H_{i}\left(
x\right) =c_{i}\right\} \subset \mathbb{C}^{6},$\ we obtain
one-parameter families which are parameterized by a curve :
\begin{equation}\label{eqn:euler}
\mathcal{D} : \theta ^{2}+c_{1}\beta ^{2}\gamma ^{2}+c_{2}\alpha
^{2}\gamma ^{2}+c_{3}\alpha ^{2}\beta ^{2}+c_{4}\alpha \beta
\gamma =0,
\end{equation}
where $\theta $\ is an arbitrary parameter and where $\alpha
=x_{4}^{\left( 0\right) },\beta =x_{5}^{\left( 0\right) },\gamma
=x_{6}^{\left( 0\right) }$parameterizes the elliptic curve
\begin{equation}\label{eqn:euler}
\mathcal{E}:\beta ^{2}=d_{1}^{2}\alpha ^{2}-1,\text{ }\gamma ^{2}=
d_{2}^{2}\alpha^{2}+1,
\end{equation}
with $d_{1},d_{2}$\ such that: $d_{1}^{2}+d_{2}^{2}+1=0.$\ The
curve $\mathcal{D}$\ is a $2$-sheeted ramified covering of the
elliptic curve $\mathcal{E}$\ $.$\ The branch points are defined
by the $16$\ zeroes of $c_{1}\beta ^{2}\gamma ^{2}+c_{2}\alpha
^{2}\gamma ^{2}+c_{3}\alpha ^{2}\beta ^{2}+c_{4}\alpha \beta
\gamma $\ on $\mathcal{E}$. The curve $\mathcal{D}$\ is unramified
at infinity and by Hurwitz's formula, the genus of $\mathcal{D}$\
is $9$. Upon putting $\zeta \equiv \alpha ^{2},$\ the curve
$\mathcal{D}$\ can also be seen as a $4-$sheeted unramified
covering of the following curve of genus $3:$
$$C:\left( \theta ^{2}+c_{1}\beta ^{2}\gamma ^{2}+
\left( c_{2}\gamma ^{2}+c_{3}\beta ^{2}\right) \zeta \right) ^{2}-
c_{4}^{2}\zeta \beta ^{2}\gamma^{2}=0.$$ Moreover, the map $\tau
:C\rightarrow C,$ $(\theta ,\zeta )\mapsto (-\theta ,\zeta ),$ is
an involution on $C$\ and the quotient $C_{0}=C/\tau $\ is an
elliptic curve defined by
$$C_{0}:\eta ^{2}=c_{4}^{2}\zeta \left( d_{1}^{2}d_{2}^{2}\zeta ^{2}+
\left( d_{1}^{2}-d_{2}^{2}\right) \zeta -1\right).$$ The curve
$C$\ is a double ramified covering of $C_{0},$ $C\rightarrow
C_{0},(\theta ,\eta ,\zeta )\mapsto (\eta ,\zeta ),$
$$
C:\left\{\begin{array}{rl} &\theta ^{2}=-c_{1}\beta ^{2}\gamma
^{2}-\left( c_{2}\gamma ^{2}+c_{3}\beta ^{2}\right) \zeta +\eta \\
&\eta ^{2}=c_{4}^{2}\zeta \left( d_{1}^{2}d_{2}^{2}\zeta
^{2}+\left( d_{1}^{2}-d_{2}^{2}\right) \zeta -1\right).
\end{array}\right.
$$
Let $(a_{1},a_{2},a_{3},b_{1},b_{2},b_{3})$\ be a canonical
homology basis of $C$\ such that $\tau \left( a_{1}\right)
=a_{3},$\ $\tau \left( b_{1}\right) =b_{3},$\ $\tau \left(
a_{2}\right) =-a_{2}$\ and $\tau \left( b_{2}\right) =-b_{2}$\ for
the involution $\tau .$ Using the Poincar\'{e} residu map, we show
that
$$\omega _{0}=\frac{d\zeta }{\eta },\text{ }\omega _{1}=
\frac{\zeta d\zeta }{\theta \eta },\text{ }\omega _{2}=
\frac{d\zeta }{\theta \eta},$$ form a basis of holomorphic
differentials on $C$\ and $\tau ^{*}\left( \omega _{0}\right)
=\omega _{0},\tau ^{*}\left( \omega _{k}\right) =-\omega _{k}$\
$\left( k=1,2\right).$ The flow evolues on an abelian surface
$\widetilde{M}_c\subseteq \mathbb{CP}^7$ of period matrix $\left(
\begin{array}{llll}
2 & 0 & a & c \\
0 & 4 & c & b
\end{array}
\right),\text{ Im}\left(
\begin{array}{ll}
a & c \\
c & b
\end{array}
\right) >0.$ Following the method (theorem 13), we obtain
\begin{Theo}
The abelian surface $\widetilde{M}_c$ can be identified as
$Prym(C/C_0)$. More precisely
$$
\overset{4}{\underset{i=1}{\bigcap }}\left\{x\in \mathbb{C}^{6},
H_{i}\left( x\right) =c_{i}\right\}=Prym(C/C_0) \backslash
\mathcal{D},$$ where $\mathcal{D}$ is a genus 9 curve (47), which
is a ramified cover of an elliptic curve $\mathcal{E}$ (48) with
16 branch points.
\end{Theo}

\section{Generalized algebraic completely integrable systems}

Some others integrable systems appear as coverings of algebraic
completely integrable systems. The manifolds invariant by the
complex flows are coverings of abelian varieties and these systems
are called algebraic completely integrable in the generalized
sense.\\
Consider the case $F_3=0,$ (see section 3.1) and the following
change of variables
$$z_{1}=q_{1}^2,\quad z_{2}=q_{2}, \quad z_{3}=p_{2},
\quad z_{4}=p_{1}{q_{1}},\quad
z_{5}=p_{1}^{2}-q_{1}^{2}q_{2}^{2}.$$ Substituting this into the
constants of motion $F_1,F_2,F_3$ leads obviously to the relations
\begin{eqnarray}
H_1&=&\frac{1}{2}p_{1}^{2}-\frac{3}{2}q_{1}^{2}q_{2}^{2}+\frac{1}{2}p_{2}^{2}
-\frac{1}{4}q_{1}^{4}-2q_{2}^{4},\\
H_2&=&p_{1}^{4}-6q_{1}^{2}q_{2}^{2}p_{1}^{2}+q_{1}^{4}q_{2}^{4}
-q_{1}^{4}p_{1}^{2}+q_{1}^{6}q_{2}^{2}+\allowbreak
4q_{1}^{3}q_{2}p_{1}p_{2}-q_{1}^{4}p_{2}^{2}+\frac{1}{4}q_{1}^{8},\nonumber
\end{eqnarray}
whereas the last constant leads to an identity. Using the
differential equations (28) combined with the transformation above
leads to the system of differential equations
\begin{eqnarray}
\ddot q_{1}&=&q_{1}\left(q_{1}^{2}+3q_{2}^{2}\right),\\
\ddot q_{1}&=&q_{2}\left( 3q_{1}^{2}+8q_{2}^{2}\right).\nonumber
\end{eqnarray}
The last equation (28) for $z_{5}$ leads to an identity. Thus, we
obtain the potential constructed by Ramani, Dorozzi and
Grammaticos [25,7]. Evidently, the functions $H_1$ and $H_2$
commute : $\left\{ H_{1},H_{2}\right\}=0.$ The system (50) is
weight-homogeneous with $q_1, q_2$ having weight 1 and $p_1, p_2$
weight 2, so that $H_1$ and $H_2$ have weight 4 and 8
respectively. When one examines all possible singularities, one
finds that it possible for the variable $q_1$ to contain square
root terms of the type $t^{1/2}$, which are strictly not allowed
by the Painlevé test (i.e., the general solutions have no movable
singularities other than poles). However, these terms are
trivially removed by introducing the variables $z_1,\ldots,z_5$
which restores the Painlevé property to the system. Let B be the
affine variety defined by
\begin{equation}\label{eqn:euler}
B=\bigcap_{k=1}^{2}\{z\in\mathbb{C}^4:H_k(z)=b_k\},\end{equation}
where $(b_{1},b_{2}) \in \mathbb{C}^{2}$.
\begin{Theo}
a) The system (50) admits Laurent solutions in $t^{1/2}$,
depending on 3 free parameters: $u, v $ and $w$. These solutions
restricted to the surface $B$(51) are parameterized by two copies
$\Gamma_{1}$
and $\Gamma_{-1}$ of the same Riemann surface of genus 16.\\
b) The invariant surface $B$(51) can be completed as a cyclic
double cover $\overline{B}$ of the abelian surface
$\widetilde{A}$, ramified along the divisor
$\mathcal{C}_{1}+\mathcal{C}_{-1}.$ The system (50) is algebraic
complete integrable in the generalized sense. Moreover,
$\overline{B}$ is smooth except at the point lying over the
singularity (of type $A_3$) of $\mathcal{C}_{1}+\mathcal{C}_{-1}$
and the resolution $\widetilde{B}$ of $\overline{B}$ is a surface
of general type with invariants : $\mathcal{X}(\widetilde{B})=1$
and $p_g(\widetilde{B})=2.$
\end{Theo}
\emph{Proof}. a) The system (50) possesses 3-dimensional family of
Laurent solutions (principal balances) depending on three free
parameters $u, v $ and $w$. There are precisely two such families,
labeled by $\varepsilon =\pm 1,$ and they are explicitly given as
follows
\begin{eqnarray}
q_{1}&=&\frac{1}{\sqrt{t}}( u-\frac{1}{4}u^{3}t+vt^{2}
-\frac{5}{128}u^{7}t^{3}+\frac{1}{8}u( \frac{3}{4}u^{3}v
-\frac{7}{256}u^{8}+\allowbreak 3\varepsilon w) t^{4}+\cdots),\nonumber\\
q_{2}&=&\frac{1}{t}( \frac{1}{2}\varepsilon
-\frac{1}{4}\varepsilon u^{2}t +\frac{1}{8}\varepsilon
u^{4}t^{2}+\frac{1}{4}\varepsilon u( \frac{1}{32}u^{5}
-3v) t^{3}+\allowbreak wt^{4}+\cdots),\\
p_{1}&=&\frac{1}{2t\sqrt{t}}( -u-\frac{1}{4}u^{3}t+3vt^{2}
-\frac{25}{128}t^{3}u^{7}+\nonumber\\
&&\qquad \qquad \qquad\frac{7}{8}u( \frac{3}{4}u^{3}v
-\frac{7}{256}u^{8}+3\varepsilon w) t^{4}+\cdots), \nonumber\\
p_{2}&=&\frac{1}{t^{2}}( -\frac{1}{2}\varepsilon
+\frac{1}{8}\varepsilon u^{4}t^{2}+\frac{1}{2}\varepsilon u(
\frac{1}{32}u^{5}-3v) t^{3}+3wt^{4}+\cdots). \nonumber
\end{eqnarray}
These formal series solutions are convergent as a consequence of
the majorant method. By substituting these series in the constants
of the motion $H_{1}=b_{1}$ and $H_{2}=b_{2},$ one eliminates the
parameter $w $ linearly, leading to an equation connecting the two
remaining parameters $u$ and $v$ :
\begin{eqnarray}
\Gamma :
&&\frac{65}{4}uv^{3}+\frac{93}{64}u^{6}v^{2}+\frac{3}{8192}\left(
-9829u^{8}+26112H_{1}\right)
u^{3}v\\&&-\frac{10299}{65536}u^{16}-\allowbreak
\frac{123}{256}H_{1}u^{8}+H_{2}+\frac{15362\,98731}{52}=0.\nonumber
\end{eqnarray}
According to Hurwitz' formula, this defines a Riemann surface
$\Gamma$ of genus 16. The Laurent solutions restricted to the
surface $B$(51) are thus parameterized by two copies $\Gamma_{-1}$
and $\Gamma_{1}$ of the same Riemann surface $\Gamma$.\\
b) The morphism $\varphi : B\longrightarrow A,\quad
(q_1,q_2,p_1,p_2)\longmapsto (z_1,z_2,z_3,z_4,z_5),$ maps the
vector field (50) into an algebraic completely integrable system
(1) in five unknowns and the affine variety $B$(51) onto the
affine part $A$(30) of an abelian variety $\widetilde{A}$ with
$\widetilde{A}\backslash A=\mathcal{C}_{1}+\mathcal{C}_{-1}$.
Observe that $\varphi$ is an unramified cover. The Riemann surface
$\Gamma$(53) play an important role in the construction of a
compactification $\overline{B}$ of $B.$ Let us denote by $G$ a
cyclic group of two elements $\{-1,1\}$ on $V_\varepsilon
^j=U_\varepsilon ^j \times \{\tau \in \mathbb{C} :
0<|\tau|<\delta\},$ where $\tau=t^{1/2}$ and $ U_\varepsilon ^j$
is an affine chart of $\Gamma_\varepsilon$ for which the Laurent
solutions (52) are defined. The action of $G$ is defined by
$(-1)\circ (u,v,\tau)=(-u,-v,-\tau)$ and is without fixed points
in $V_\varepsilon ^j.$ So we can identify the quotient
$V_\varepsilon ^j / G$ with the image of the smooth map
$h_\varepsilon ^j :V_\varepsilon ^j\rightarrow B$ defined by the
expansions (10). We have $(-1,1).(u,v,\tau)=(-u,-v,\tau)$ and
$(1,-1).(u,v,\tau)=(u,v,-\tau),$ i.e., $G\times G$ acts separately
on each coordinate. Thus, identifying $V_\varepsilon ^j/G^2$ with
the image of $\varphi\circ h_\varepsilon ^j$ in $A.$ Note that
$B_\varepsilon ^j=V_\varepsilon ^j/G$ is smooth (except for a
finite number of points) and the coherence of the $B_\varepsilon
^j$ follows from the coherence of $V_\varepsilon ^j$ and the
action of $G.$ Now by taking $B$ and by gluing on various
varieties $B_\varepsilon ^j\backslash \{\mbox{some points}\},$ we
obtain a smooth complex manifold $\widehat{B}$ which is a double
cover of the abelian variety $\widetilde{A}$ (constructed in
proposition 2.3) ramified along $\mathcal{C}_1+\mathcal{C}_{-1},$
and therefore can be completed to an algebraic cyclic cover of
$\widetilde{A}.$ To see what happens to the missing points, we
must investigate the image of $\Gamma \times\{0\}$ in $\cup
B_\varepsilon ^j.$ The quotient $\Gamma \times\{0\}/G$ is
birationally equivalent to the Riemann surface $\Upsilon$ of genus
7 :
$$\Upsilon : \frac{65}{4}y^{3}+\frac{93}{64}x^{3}y^{2}+\frac{3}{8192}
\left( -9829x^{4}+26112b_{1}\right) x^{2}y\\$$$$+ x\left(
-\frac{10299}{65536}x^{8}-\allowbreak
\frac{123}{256}b_{1}x^{4}+b_{2}+
\frac{15362\,98731}{52}\right)=0,$$ where $y=uv, x=u^2.$ The
Riemann surface $\Upsilon$ is birationally equivalent to
$\mathcal{C}.$ The only points of $\Upsilon$ fixed under
$(u,v)\mapsto (-u,-v)$ are the points at $\infty,$ which
correspond to the ramification points of the map $\Gamma
\times\{0\}\overset{2-1}{\rightarrow }\Upsilon :(u,v)\mapsto(x,y)$
and coincides with the points at $\infty$ of the Riemann surface
$\mathcal{C}.$ Then the variety $\widehat{B}$ constructed above is
birationally equivalent to the compactification $\overline{B}$ of
the generic invariant surface $B.$ So $\overline{B}$ is a cyclic
double cover of the abelian surface $\widetilde{A}$ ramified along
the divisor $\mathcal{C}_1+\mathcal{C}_{-1},$ where
$\mathcal{C}_1$ and $\mathcal{C}_{-1}$ have two points in commune
at which they are tangent to each other. It follows that The
system (8) is algebraic complete integrable in the generalized
sense. Moreover, $\overline{B}$ is smooth except at the point
lying over the singularity (of type $A_3$) of
$\mathcal{C}_1+\mathcal{C}_{-1}.$ In term of an appropriate local
holomorphic coordinate system $(X,Y,Z),$ the local analytic
equation about this singularity is $X^4+Y^2+Z^2=0.$ Now, let
$\widetilde{B}$ be the resolution of singularities of
$\overline{B},$ $\mathcal{X}(\widetilde{B})$ be the Euler
characteristic of $\widetilde{B}$ and $p_g(\widetilde{B})$ the
geometric genus of $\widetilde{B}.$ Then $\widetilde{B}$ is a
surface of general type with invariants :
$\mathcal{X}(\widetilde{B})=1$ and $p_g(\widetilde{B})=2.$ This
concludes the proof of the theorem.

\begin{Rem} The asymptotic solution (52) can be read off from (31) and the
change of variable : $q_1=\sqrt{z_1}, q_2=z_2, p_1=z_4/q_1,
p_2=z_3.$ The function $z_1$ has a simple pole along the divisor
$\mathcal{C}_1+\mathcal{C}_{-1}$ and a double zero along a Riemann
surface of genus 7 defining a double cover of $\widetilde{A}$
ramified along $\mathcal{C}_1+\mathcal{C}_{-1}.$
\end{Rem}

\end{document}